\documentclass[11pt,reqno]{amsart}
\usepackage{amssymb}
\usepackage{amsmath}
 \usepackage[usenames,dvipsnames]{color}
\usepackage[kerning=true]{microtype}
\usepackage[page]{appendix}
\usepackage[all,cmtip]{xy}
\usepackage{pgf,tikz}
 \usetikzlibrary{arrows,shapes}
  \usetikzlibrary{calc}

\usepackage{fullpage}
\usepackage{amsthm}
\usepackage{bm}
\usepackage{MnSymbol}
\usepackage[T1]{fontenc}
\usepackage[applemac]{inputenc}
\usepackage[colorlinks=true, citecolor=cyan,urlcolor=blue,%linktocpage=true,
pagebackref=true,hypertexnames=true]{hyperref} 
\usepackage{genyoungtabtikz}

\YFrench

\newcommand\yred{\Yfillcolour{red}}
\newcommand\ylw{\Yfillcolour{yellow}}
\newcommand\ygreen{\Yfillcolour{green!80}}
\newcommand\ywhite{\Yfillcolour{white}}
  \Yboxdim{15pt}
\Yfillcolour{green!80}

\definecolor{azure}{rgb}{0.,0.5,1.0} 
\definecolor{green}{rgb}{0.,0.7,0.} 

\def\bleu{\textcolor{blue}}

\newcommand{\auteur}[1]{{\sc #1}}
\newcommand{\titreref}[1]{{\em #1}}
\newcommand{\vol}[1]{{\bf #1}}
\newcommand{\MathReview}[1]{\href{http://www.ams.org/mathscinet-getitem?mr=#1}{MR#1}.}
\renewcommand{\MathReview}[1]{ }

\def\charac{\raise 2pt\hbox{\large$\chi$}}
\newcommand{\Cjn}[2]{C_{#1}^{\,#2}}
\DeclareMathOperator*{\leftmost}{\leftarrow}
\DeclareMathOperator*{\rightmost}{\rightarrow}

\newcommand{\GL}{\mathrm{GL}}
\newcommand{\Id}{\mathrm{Id}}

\newcommand{\MacH}{\widetilde{H}}
\newcommand{\M}{\mathcal{M}}

\newcommand{\WN}{\mathcal{W}}

\newcommand{\myatop}[2]{\genfrac{}{}{0pt}{}{#1}{#2}}
\newcommand{\N}{\mathbb{N}}
\newcommand{\qbinom}[2]{\genfrac{[}{]}{0pt}{}{#1}{#2}_q}
\newcommand{\tqbinom}[2]{\textstyle \genfrac{[}{]}{0pt}{}{#1}{#2}_q}

\newcommand{\Rational}{\mathbb{Q}}
\newcommand{\R}{\mathcal{R}}
\newcommand{\scalar}[2]{{\langle#1,#2 \rangle}}

\renewcommand{\S}{\mathbb{S}}
\newcommand{\sgn}{\mathrm{sgn}}
\newcommand{\pref}[1]{{\rm (\ref{#1})}}

\newcommand{\define}[1]{\bleu{\bf{#1}}}
\newcommand\zero{{\scriptscriptstyle 0}}

\def\flechegauche#1{\xymatrixcolsep{#1pc}\xymatrix{&\ar[l] &}}
\def\flechedroite#1{\xymatrixcolsep{#1pc}\xymatrix{\ar[r] &}}
\def\flechebas#1{\xymatrix@R=#1pc{& \ar[d] \\  &}}
\def\flechehaut#1{\xymatrix@R=#1pc{\\ \ar[u] &}}

\def\epi#1#2{\xymatrixcolsep{1.2pc}\xymatrix{#1\,\ar@{->>}[r]&\,#2}}

\def\mono#1#2{\xymatrixcolsep{1pc}\xymatrix{#1\ \ar@{^{(}->}[r] &\,#2}}

\newtheorem{lemma}{\bleu{Lemma}}

\newtheorem{prop}{\bleu{Proposition}}
\newtheorem{definition}{\bleu{Definition}}

\numberwithin{equation}{section}
\numberwithin{lemma}{section}
\numberwithin{rmk}{section}
\numberwithin{cor}{section}

\title{A Survey  of $q$-Whittaker polynomials}
\author{F.~Bergeron}

\address{\href{http://bergeron.math.uqam.ca}{D\'epartement de Math\'ematiques, Lacim, UQAM.}}
\email{\href{mailto:bergeron.francois@uqam.ca}{bergeron.francois@uqam.ca}}
\urladdr{bergeron.math.uqam.ca}

\thanks{During this research, the author was supported by a NSERC grant.}
\date{\bleu{\bf \today}}

\begin{document}

\begin{abstract} 
Exploiting the fact that the $q$-Whittaker polynomials arise as a specialization of the (modified) Macdonald polynomials, we derive some of their basic properties, and explore interesting identities that they satisfy. We also show how they arise as graded Frobenius characteristics of $\S_n$-modules, and give a combinatorial approach to associated Pieri formulas.
\end{abstract}

\maketitle
	{ \setcounter{tocdepth}{1}\parskip=0pt\footnotesize \tableofcontents}
\parskip=5pt  
\parindent=15pt

%%%%%%%%%%%%%%%%%%
%              Section                          %
%%%%%%%%%%%%%%%%%%
\section{\bleu{Introduction}}
Introduced under this name in~\cite{Gerasimov}, the $q$-Whittaker polynomials arise as a specialization of the Macdonald polynomials (see~\cite{macdonald_lotha} and appendix), and occur in the context of the study of ``common eigenfunction of a set of commuting $q$-deformed Toda chain Hamiltonians''. However, they were studied much earlier under other guises, since they are directly related to Hall-Littlewood polynomials, and are explicitly discussed in Macdonald~\cite{macdonald}. The purpose of these notes is to survey some of their basic properties, and to discuss interesting questions/identities concerning them. Moreover, we describe how they occur as graded  Frobenius characteristics of explicit $\S_n$-modules.

%%%%%%%%%%%%%%%%%%%%%%%%
\section{\bleu{Combinatorics of diagram}}\label{section_notations}
For symmetric functions, and related notions, we mostly follow Macdonald's notations, except for diagrams for which we consider more natural to use the French convention. As usual, a length $k$ partition $\mu$ of $n$ is a list of integers $\mu=\mu_1\mu_2\cdots \mu_k$, with $\mu_1\geq \mu_2\geq \ldots \geq \mu_k\geq 1$ and $n=|\mu|:=\mu_1+\mu_2+\ldots+\mu_k$. We sometimes add a part $\mu_0=\infty$ to $\mu$, and infinitely many $0$-parts to $\mu$, so that $\mu$ may take (as needed) one of the forms
   $$\mu=(\mu_1,\ldots,\mu_k) = (\infty,\mu_1,\ldots,\mu_k)=(\infty,\mu_1,\ldots,\mu_k,0,0,\ldots).$$
We call \define{cell} any element of $\N\times \N$, and say that $c=(i,j)$ is a cell of $\mu$ whenever $0\leq i+1\leq \mu_{j+1}$. We then write $c\in\mu$. Cells of $\mu$ correspond to Cartesian coordinates of southwest corners of the $1\times 1$ boxes lying in the \define{diagram} of $\mu$ (in French notation). Thus $(0,0)$ is the southwest-most cell in (the diagram of) $\mu$. A cell $c=(i,j)$ in $\mu$ lies on the \define{row} $j+1$ of the diagram of $\mu$. 

\begin{figure}[ht]
\begin{center}
\setlength{\unitlength}{1em}
\begin{picture}(30,10)
\put(0,0){\yng(14,9,9,9,8,6,6)}
\put(2.75,2.75){\gyoung(!\yred;!\ylw;;;;;;,;,;,;,;)}
\put(3.1,2.9){$\bm{c}$}
\put(3.8,2.8){$\flechegauche{2.5}$}
	\put(7.5,2.9){\tiny$\bm{a(c)}$}
	\put(8.8,2.8){$\flechedroite{2.5}$}
\put(0.3,5.9){$\flechebas{1.6}$}
\put(2.8,6.1){\tiny$\bm{\ell(c)}$}
\put(3.1,9.3){$\flechehaut{2}$}
\put(15.3,0){\put(3.8,0){$\flechegauche{2.5}$}
	\put(7.6,0){\small$\sigma_\mu(0)$}
	\put(10,0){$\flechedroite{2.5}_\infty$}}
\put(-4.2,10){\put(3.8,0){$\flechegauche{2}$}
	\put(7,0){\small $\sigma_\mu(k)$}
	\put(9,0){$\flechedroite{2.4}$}}
\put(8.2,1.5){\put(3.8,0){$\flechegauche{1.5}$}
	\put(6.5,0){\small $\sigma_\mu(1)$}
	\put(8.5,0){$\flechedroite{1.8}$}}
\end{picture}
\end{center}
\qquad  \vskip-15pt
\caption{Some of the step sizes of a partition; and the arm and leg of a cell $c$ in $\mu$.}
\label{Fig_part} 
\end{figure}
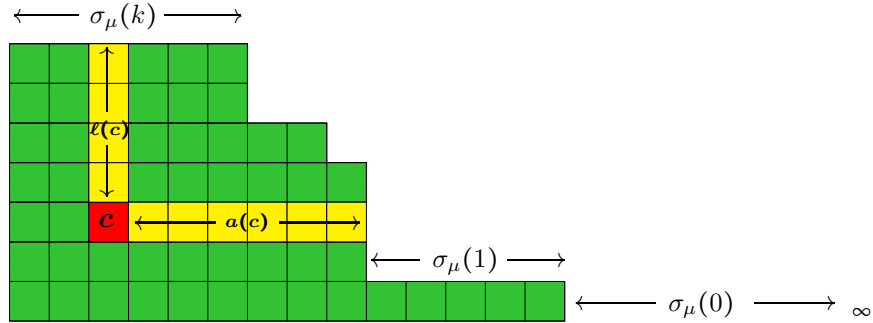

The \define{cell enumerator} polynomial $B_\mu(q,t)$ is defined as $B_\mu(q,t):=\sum_{(i,j)\in \mu} q^it^j$.
  For example,
  	$$B_{52}(q,t)= 1+q+q^2+q^3+q^4+t+qt. $$
The $i^{\rm th}$-\define{step size} of a partition $\mu$, denoted by $\sigma_\mu(i)$, is the part difference $\mu_i-\mu_{i+1}$. In particular,  the $0^{\rm th}$-\define{step size} is always equal to $\sigma_\mu(0):=\infty$. Clearly, step size are all equal to $0$ for $0$-parts.
Clearly, some steps may be equal to $0$, and the last non-zero step is equal to the smallest part of $\mu$.
The \define{step sequence} of $\mu$ is $\sigma(\mu):=(\sigma_\mu(0),\sigma_\mu(1),\sigma_\mu(2),\ldots,\sigma_\mu(k))$, omitting step sizes of $0$-parts. For example, the step sequence of the partition in Figure~\ref{Fig_part} is $(\infty,5,0,0,1,2,0,6)$. 
For a given partition $\lambda$, we extend the usual notions of \define{arm length} and \define{leg length}   to the whole $\N\times \N$-plane,  setting
 $$a_\lambda(i,j):=\begin{cases}
      \lambda_{j+1}-(i+1)& \text{if}\ (i,j)\in\lambda, \\
     i-\lambda_{j+1} & \text{otherwise},
\end{cases}$$
 and  $\ell_\lambda(i,j):=a_{\lambda'}(j,i)$.  
 Observe that there are two kinds of cells having arm length (resp. leg length) equal to $0$. Those are the cells that lie either immediately to the left or immediately the right of the right ``boundary'' of $\mu$. In the figure below, a green-colored partition is illustrated, with each cell (inside and outside of $\mu$) marked by arm length. 
\begin{center}
\setlength{\unitlength}{1em}
 \begin{picture}(9,6)
 \Yfillcolour{white}
 \put(0,0){\gyoung(;;;;;0;{1};{2},;;;;;0;{1};{2},;;;;0;{1};{2};{3},0;{1};{2};{3};{4};{5};{6})}
 \Yfillcolour{green!80}
 \put(0,0){\gyoung(;3;2;1;0,;3;2;1;0,;2;1;0)}
 \end{picture}
\end{center}
When possible, we write more simply $a(c)$ for $a_\lambda(c)$ and $\ell(c)$ for $\ell_\lambda(c)$.
When needed, we will also consider that we have cells $(i,-1)$ in the part infinite part $\lambda_0$ of $\lambda$, for all $i\geq 0$. There leg length are set to be 
$$\ell_\lambda(i,-1):=\begin{cases}
      \lambda'_{i+1}& \text{if}\ (i,0)\in\lambda, \\
      0 & \text{otherwise}.
\end{cases}$$
With these definitions, an \define{inner corner} (resp. \define{outer corner}) of $\mu$ is a cell such that $a(c)=\ell(c)=0$ and $c\in\mu$  (resp. $c\notin\mu$).

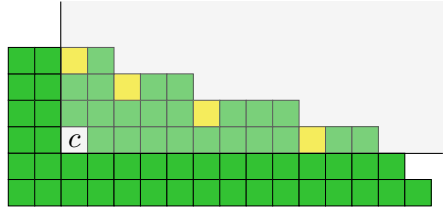
\begin{figure}[ht]
\begin{center}
\begin{tikzpicture}
  \Yboxdim{10pt}
\node at (.535,-.165) {\gyoung(;;;;;;;;;;;;;;;;,;;;;;;;;;;;;;;;,;;!\ywhite;!\ygreen;;;;;;;;!\ylw;!\ygreen;;,;;;;;;;!\ylw;!\ygreen;;;,;;;;!\ylw;!\ygreen;;,;;!\ylw;!\ygreen;)};
\path [fill=gray!20,opacity=.4] (-1.58,-.52) rectangle (3.5,1.5);
\node at (-1.4,-.35) {$c$};
\draw (-1.58,-.52) -- (-1.58,1.5);
\draw (-1.58,-.52) -- (3.5,-.52);
\end{tikzpicture}
\end{center}\caption{The cells $c'$ (marked in yellow) such that $c\leftmost_\lambda c'$.}\label{Figure_WDc}
\end{figure}

 It will be handy, for a given cell $c$, to consider all the cells $c'$ in $\lambda$ which satisfy the following conditions:
\begin{itemize}
\item $c'\in\lambda$,
\item  $\ell(c')=0$, 
 \item $c'$ lie to the north-east of $c$, and  
 \item $c'$ is the leftmost such cell on its row. 
 \end{itemize}
 We write $c\leftmost_{\lambda} c'$ when this is the case.
For instance, for a  given cell  $c$, the cells such that $c\leftmost_\lambda c'$ are marked in yellow\footnote{All lying in the shaded region which corresponds to the condition $c\preceq \gamma$.} in Figure~\ref{Figure_WDc}. 
\noindent Naturally, the set $\{c'\,|\, c\leftmost_\lambda c'\}$ is empty when $c$ does not lie in $\lambda$. 
Observe that, when $c=(i,j)$ is such that $(0,0)\leftmost_\lambda c$, we have 
\begin{equation}\label{step_arm}
   \sigma_\lambda(j+1)=a(c)+1,\qquad{\rm hence}\qquad
          [\sigma_\lambda(j+1)]_q=1+q+\ldots+q^{a(c)}. 
  \end{equation}
  \begin{figure}[ht]
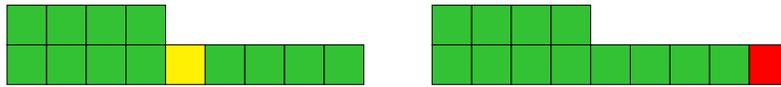

 \begin{center}
 \gyoung(;;;;!\ylw;!\ygreen;;;;,;;;;) \qquad 
 \gyoung(;;;;;;;;!\yred;!\ygreen,;;;;)
  \end{center}\caption{An internal corner (yellow) vs an inner corner (red), for the partition $(9,4)$.}
 \end{figure}
These cells $c$ are said to be \define{internal corners} of $\lambda$. For each internal corner $c=(i,j)$ of $\lambda$ there is an associated inner corner $c'=(i',j)$ of $\lambda$ on the same row, for which $i'=i+a(c)$.

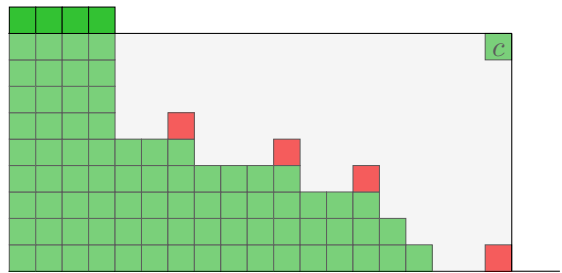
\begin{figure}[ht]
\begin{center}
\begin{tikzpicture}
  \Yboxdim{10pt}
 \draw (-2.83,-1.922) -- (4.5,-1.922);
\draw (-2.83,1.242) -- (3.86,1.242);
\draw (3.86,-1.922) -- (3.86,1.242);
\node at (.535,-.165) {\gyoung(;;;;;;;;;;;;;;;;::!\yred;,!\ygreen;;;;;;;;;;;;;;;,;;;;;;;;;;;;;;,;;;;;;;;;;;::!\yred;!\ygreen,;;;;;;;:::!\yred;!\ygreen,;;;;::!\yred;!\ygreen,;;;;,;;;;,;;;;::::::::::::::!\ygreen;c,!\ygreen;;;;)};
\path [fill=gray!20,opacity=.4] (-2.83,-1.922) rectangle (3.85,1.23);
\end{tikzpicture}\caption{The cells $c'$ (marked in red) such that $c'\rightmost_\lambda c$.}\label{Figure_WUc}
\end{center}
\end{figure} 
Similarly for $c$ in $\N^*\times \N^*$, with $\N^*=\N\cup\{\infty\}$, we consider the cells $c'$ characterized by the properties
\begin{itemize}
 \item  $c'\not\in\lambda$,  
 \item  $\ell(c')=0$,  
 \item $c'$ lie to the south-west of $c$, and
 \item $c'$ is the rightmost such cell on its row. 
 \end{itemize}
We then write $c'\rightmost_\lambda c$, and include the case $c'=(\infty,0)$ if $c=(\infty,j)$, considering that $\ell(\infty,0)=0$. The set $\{c'\,|\, c'\rightmost_\lambda c\}$ is clearly empty when $c\in\lambda$. This is illustrated in Figure~\ref{Figure_WUc}.
For the cells $c=(i,j)$ such that $c\rightmost_\lambda(\infty,\infty)$, we 
also have~\pref{step_arm},
 considering that $a(\infty,0)=\infty$. These cells $c$ are said to be \define{external corners} of $\lambda$. For each external corner $c=(i,j)$ of $\lambda$ there is an associated outer corner $c'=(i',j)$ of $\lambda$ on the same row (immediately to left of $\lambda$), for which $i=i'+a(c)$.
 \begin{figure}[ht]
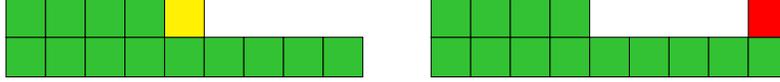

 \begin{center}
 \gyoung(;;;;;;;;;!\ygreen,;;;;!\ylw;) \qquad 
 \gyoung(;;;;;;;;;!\ygreen,;;;;::::!\yred;)
  \end{center}\caption{An outer corner (yellow) vs an external corner (red), for the partition $(9,4)$.}
 \end{figure}

\subsection*{About $q$-Analogs}
\begin{figure}[ht]
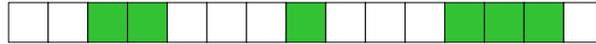

\begin{center}
\gyoung(!\ywhite;;!\ygreen;;!\ywhite;;;!\ygreen;!\ywhite;;;!\ygreen;;;!\ywhite;)
\end{center}\caption{A $k$-subset, with cells in green.}\label{Fig_ksubdiag}
\end{figure}

Recall that the usual $q$-\define{analog} of $n$ is $[n]_q:=1+q+\ldots +q^{n-1}$. It is natural to set 
	$$[\infty]_q:=\frac{1}{1-q}.$$
This makes it natural to extend the usual $q$-\define{binomial}:
  $$\qbinom{n}{k} :=\frac{[n]_q\cdots [n-k+1]_q}{[k]_q\cdots [1]_q},\qquad n\in\N,$$
 to the infinite context, by setting
  $$\qbinom{\infty}{k}:=\prod_{i=1}^{r}\frac{1}{1-q^i}.$$
To better see this, as well as connections to our discussion, it may be worth recalling that $\qbinom{n}{k}$ $q$-counts 
$k$-subsets of the $n$-set $\{0,1,\ldots,n-1\}$, which may clearly be represented as diagrams such as illustrated in Figure~\ref{Fig_ksubdiag}
This is to say that
\begin{equation}
   \qbinom{n}{k}=\sum_{|\bm{d}|=k}  q^{\Sigma(\bm{d})-\binom{k}{2}},\qquad {\rm with}\qquad
     \Sigma(\bm{d}):=\sum_{i\in \bm{d}} i.
\end{equation}
Thus the $q$-weight degree of each cell $(i,j)$ is equal to $i$ minus the number of cells to its left. It follows that the limit as $n$ tends to $\infty$ corresponds bijectively to partitions at most $k$ parts, via the usual bijection that sends $0\leq a_0<a_1<\ldots < a_{k-1}$ to the partition having parts equal to $a_i-i$. One may also think\footnote{See for instance~\cite{borodin_corwin,williams} for related notions of statistical mechanics.} of this as a global slide ``displacement'' measure. For this, one starts with all cell adjacent and justified to the right inside an overall set: $\{n-k,n-k+1,\ldots,n-1\}$. Cells are then allowed to slide to the left in all possible way, one unit at a time, except that cells are forbidden to cross over one another. The $q$-weight degree then measure how many single slides are needed to go from the original configuration to a given one. 

\begin{figure}[ht]
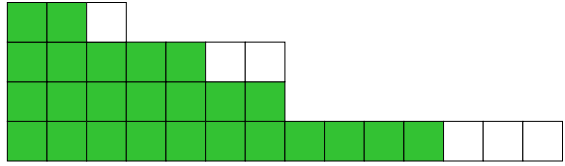

 \begin{center}
 \gyoung(;;;;;;;;;;;!\ywhite;;;!\ygreen,;;;;;;;,;;;;;!\ywhite;;!\ygreen,;;!\ywhite;) \qquad 
  \end{center}\caption{The diagram $(10,7,5,2)$ inside $(13,7,7,3) $.}\label{Fig_slide_diag}
 \end{figure}

The following $q$-polynomial arises when one extends the above sliding process to diagrams such as the one illustrated in Figure~\ref{Fig_slide_diag}
 \begin{equation}\label{c_mu_nu}
 c_{\mu\lambda}(q):=\prod_{i\geq 1} \qbinom{\sigma_\lambda(i)}{\lambda_i-\mu_i}.
 \end{equation}
 For sure this is not zero only of $\mu$ is contained in $\lambda$.
We assume given two partitions $\mu$ and $\lambda$, such that the skew shape $\lambda/\mu$ is an \define {horizontal strip}. Thus, for all $i$, $\lambda_i-\mu_i\leq \sigma_\lambda(i)$. Inside the shape $\lambda$ one considers the cells of $\mu$. On each row, the last $\sigma_\lambda(i)-(\lambda_i-\mu_i)$ of these cells of $\mu$ are slided as far right as possible inside the overall shape $\lambda$, with no overlap. The resulting \define{shifted diagram} is denoted $\mu\rightarrow \lambda$. 
\begin{figure}[ht]
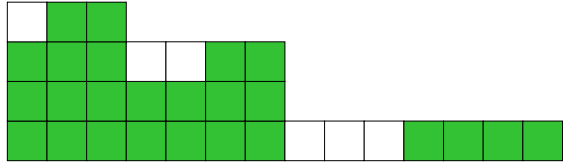

 \begin{center}
 \gyoung(;;;;;;;!\ywhite;;;!\ygreen;;;;,;;;;;;;,;;;!\ywhite;;!\ygreen;;,!\ywhite;!\ygreen;;) \qquad 
  \end{center}\caption{The shifted diagram $(10,7,5,2)\rightarrow(13,7,7,3) $.}
 \end{figure}
 \begin{lemma}
 For all partitions $\mu$ and $\lambda$, we have
     \begin{equation}
         c_{\mu\lambda}(q)=\sum_{\bm{d}} q^{{\rm slide}_{\mu\rightarrow \lambda}(\bm{d})},
       \end{equation}
   where the sum is over all diagrams that may be obtained from $\mu\rightarrow \lambda$ by sliding cells to the left, with no overlap. The ${\rm slide}_{\mu\rightarrow \lambda}(\bm{d})$ degree is the number of unit slides needed to go from $\mu\rightarrow \lambda$ to ${\bm d}$. 
 \end{lemma}

In the description of the dual Whittaker polynomial (see~\pref{definition_W_dual}), we will make use of the following weight:
	\begin{equation}\label{top_cells_q_arms}
	     v_\mu(q):= \prod_{\ell(c)=0} (1-q^{a(c)+1}) =(1-q)^{\mu_1}\prod_{i=1}^k \sigma_\mu(i)!_q,
	\end{equation}
for $c\in\mu$.
One may check that (using plethystic notation) 
    \begin{equation}
       h_n\Big[\frac{1}{1-q}\Big] =h_n(1,q,q^2,\ldots)= \prod_{k=1}^n \frac{1}{1-q^k},\qquad {\rm hence}\qquad   h_{\sigma(\mu)}\Big[\frac{1}{1-q}\Big]=\frac{1}{v_\mu(q)}.
    \end{equation}

%%%%%%%%%%%%%%%%%%%%%%

\section{\bleu{\texorpdfstring{$q$}{q}-Whittaker polynomials}}
Equivalent descriptions of the $q$-Whittaker polynomials, denoted by $W_\mu(q;\bm{z})$, include the following\footnote{See appendix for more notations and definitions. In paper~\cite{borodin}, one may find connections to physics.}:
 \begin{align*}
    W_\mu(q;\bm{z})=\MacH_\mu(q,t;\bm{z})\big|_{t^{\eta(\mu)}}
       =H_\mu(q,0;\bm{z})
       =P_\mu(q,0;\bm{z})
       =\omega\, q^{\eta(\mu')} \MacH_\mu(1/q,0;\bm{z})
       =\omega\, Q'_{\mu'}(q;\bm{z})
 \end{align*}
where $\omega$ is the classical involution on symmetric functions, that send $s_\mu$ to $s_{\mu'}$.
 The polynomials $\mathcal{H}_\mu(q;\bm{z}) = \MacH_\mu(q,0;\bm{z})$, are sometimes called the \define{Hall-Littlewood}\footnote{But there are also other slightly different conventions.} polynomials. They appear as Frobenius transforms of the graded character of the singular co-homology ring of Springer fibers (corresponding to nilpotent matrices of Jordan type $\mu$, see~\cite{de_concini_procesi,garsia_procesi,griffin,Tymoczko}). Observe that, in the special case $\mu=(n)$ we have $W_{n} =H_n=\MacH_{n}$. In general, we have the $t$-expansion:
 \begin{equation}
    \MacH_{\mu}(q,t;\bm{z}) = \mathcal{H}_{\mu}(q;\bm{z})+\ldots + t^{\eta(\mu)}\,W_{\mu}(q;\bm{z})
 \end{equation}
Many known identities regarding the above polynomials may thus clearly be translated into statements  for the $W_\mu$'s. 

The above equalities, compare  $W_{31}= s_{31}+ qs_{22} + (q^{2} + q)s_{211} + q^{3}s_{1111}$, with
\begin{align*}
&\MacH_{31}(q,t;\bm{z})= s_{4}   
	+ (q^{2} + q + t)\,s_{31}
	+ (q^{2} + q\, t)\,s_{22}
	+ (q^{3} + q^{2} t + q\, t)\,s_{211} 
	+ q^{3} t\,s_{1111},\\
&H_{31}(q,t;\bm{z})= t\,s_{4} 
	+ (q^{2} t + q\, t + 1)\,s_{31} 
	+ (q^{2} t + q)\,s_{22} 
	+ (q^{3} t + q^{2} + q)\,s_{211} 
	+q^{3}s_{1111},\\
&P_{31}(q,t;\bm{z})=s_{31}
	+ \frac{q - t}{1-q\, t }\,s_{22}
	+ \frac{(1+q)(q-t)(1-q\,t^2)}{(1-q t)(1-q^2t^2)  }\,s_{211} 
	+\frac{(1+t)(q-t)(q^2-t)}{(1-q t)(1-q^2t^2)  }\,s_{1111} \\
&\mathcal{H}_{31}(q;\bm{z})= s_{4}   
	+ (q^{2} + q)\,s_{31}
	+ q^{2}\,s_{22}
	+ q^{3}\,s_{211},\\
&Q'_{211}(q;\bm{z})=s_{211} + q\,s_{22} + (q^{2} + q)\,s_{31} + q^{3}s_{4}.
\end{align*}
In particular,
\begin{align*}
	&\MacH_{31}(q,t;\bm{z}) = \mathcal{H}_{31}(q;\bm{z})+ t\,W_{31}(q;\bm{z}),\\
	&H_{31}(q,t;\bm{z}) =t\, \mathcal{H}_{31}(q;\bm{z})+ W_{31}(q;\bm{z}),\\
	&W_{31}(q;\bm{z})= \omega\, q^3 \mathcal{H}_{31}(1/q;\bm{z}),\quad {\rm and}\\
	&W_{31}(q;\bm{z})= \omega\,  Q'_{211}(1/q;\bm{z}).
\end{align*}
The $W_\mu$'s are Schur positive, which is to say that they expand with coefficients in $\N[q]$ in the Schur basis. One has the following specializations/formulas
\begin{equation}\label{W_coefficients}
W_\mu(0;\bm{z})=s_\mu(\bm{z}),\qquad \langle W_\mu(q;\bm{z}),e_n(\bm{z})\rangle =q^{\eta(\mu')}, \qquad {\rm and}\qquad W_\mu(1;\bm{z})=e_{\mu'}(\bm{z}).
\end{equation}
The \define{dual Whittaker polynomials} are defined as
\begin{equation}\label{definition_W_dual}
    \widehat{W}_\mu:=\frac{1}{v_\mu(q,t)} W_\mu,
 \end{equation}
 with $v_\mu(q,t)$ specified in~\pref{top_cells_q_arms}.
Observe that 
\begin{align}
&\widehat{W}_{n}(q;\bm{z}) =h_n\Big[\frac{\bm{z}}{1-q} \Big],\qquad {\rm and} \label{formule_W_n}\\
& \widehat{W}_{1^n}(q;\bm{z}) =\frac{e_n(\bm{z})}{1-q}.
\end{align}
Also, for all partitions $\mu$ of $n$,  we have the formulas
\begin{align}
    W_\mu[q;1-u] 
     &=(-u)^{n}q^{\eta(\mu)}\prod_{\myatop{(i,j)\in\mu}{\ell(i,j)=0}}\big( 1-{1}/{(u\,q^i)}\big)\label{formula_hook}\\
    (\omega\,W_\mu)[q;1-u]
     &=q^{\eta(\mu)}\prod_{\myatop{(i,j)\in\mu}{\ell(i,j)=0}}\big( 1-{u}/{q^i}\big)\label{formula_hook2}
 \end{align}
From ~\pref{formula_hook}, we may easily calculate the coefficient of hook indexed $s_{(n-k,1^k)}$ in $W_\mu$, since we have the evaluation
 $$ s_\mu[1-u]=\begin{cases}
     (-u)^k(1-u) & \text{if}\ \mu=(n-k,1^k), \\
      0& \text{otherwise}.
\end{cases}$$

%%%%%%%%%%%%%%%%%%%%%%%%%%%%%%%%%%%%%%%
\subsection*{Combinatorial formula}
The following formula (See~\cite{macdonald} for the definition of the Kostka-Foulkes polynomials $K_{\lambda\mu}(q)$ in terms of charge) gives the Schur expansion of the Whittaker polynomials.
\begin{equation}\label{Comb_Formula}
	W_\mu(q;\bm{z}) 
	= \sum_{\lambda\vdash n} K_{ \lambda'\mu'}(q)\, s_{\lambda}(\bm{z})
\end{equation}
For example, we have the matrix
$$(K_{ \lambda'\mu'}(q))=\left(\begin{array}{rrrrr}
1 & q^{3} + q^{2} + q & q^{4} + q^{2} & q^{5} + q^{4} + q^{3} & q^{6} \\
0 & 1 & q & q^{2} + q & q^{3} \\
0 & 0 & 1 & q & q^{2} \\
0 & 0 & 0 & 1 & q \\
0 & 0 & 0 & 0 & 1
\end{array}\right)$$
with rows and columns indexed by partitions in lexicographic order: $4,31,22,211,1111$. 
As it happens, $W_\mu(q;\bm{z})$ is the Frobenius transform of the character of a graded module (see section~\ref{sectionmodule}). Its graded dimension, or Hilbert series, which is here denoted by $W_\mu(q)$,  may be obtained by taking the scalar product
   \begin{equation}\label{Schur_Expansion}
          W_\mu(q)=\scalar{W_\mu(q;\bm{z})}{p_1^n(\bm{z})},
    \end{equation}
 since this is the value of the (graded) character at the identity. From~\pref{Schur_Expansion}, we get that
\begin{equation}\label{W_Kostka}
   W_\mu(q)= \sum_{\lambda\vdash n} K_{ \lambda'\mu'}(q)\, f^\lambda,
\end{equation}
hence
   $$W_\mu(0)=f^\mu\quad \qquad {\rm and} \qquad W_\mu(1)=\binom{n}{\mu'},$$
  with the last term expressed as a multinomial coefficient.

  One has the formulas
    \begin{equation}\label{formula_g_mu}
        W_\mu(q)= g_\mu(q)\,\textstyle \prod_i \sigma_\mu(i)!_q, \qquad {\rm and}\qquad
       \displaystyle \widehat{W}_\mu(q)= \frac{g_\mu(q)}{(1-q)^{\mu_1}}.
    \end{equation}
 where  $g_\mu$ is some positive coefficient (monic) polynomial whose constant term is $f^\mu$. 
 
%%%%%%%%%%%%%%%%%%
%              Section                          %
%%%%%%%%%%%%%%%%%% 
\section{\bleu{Some explicit values}}
To get a better feeling of the behavior of Whittaker polynomials, here are a few of them, expanded in the Schur basis. 

{\small
$\begin{aligned}
&W_{2}=s_{2}+qs_{11}  ,\\
&W_{11}=s_{11};\\[6pt]
&W_{3}=s_{3}+ (q^{2} + q)s_{21}+ q^{3}s_{111}  ,\\
&W_{21}= s_{21}+qs_{111} ,\\
&W_{111}=s_{111};\\[6pt]
&W_{4}= s_{4} + (q^{3} + q^{2} + q)s_{31}+ (q^{4} + q^{2})s_{22}+ (q^{5} + q^{4} + q^{3})s_{211}   +q^{6}s_{1111},\\
&W_{31}= s_{31}+ qs_{22} + (q^{2} + q)s_{211} + q^{3}s_{1111}\\
&W_{22}=s_{22}+ qs_{211} + q^{2}s_{1111} ,\\
&W_{211}= s_{211}+qs_{1111} ,\\
&W_{1111}=s_{1111};\end{aligned}$

$\begin{aligned}
&W_{5}=  s_{5}+ (q^{4} + q^{3} + q^{2} + q)s_{41}+ (q^{6} + q^{5} + q^{4} + q^{3} + q^{2})s_{32}+ (q^{7} + q^{6} + 2 q^{5} + q^{4} + q^{3})s_{311} \\
&\qquad\qquad+ (q^{8} + q^{7} + q^{6} + q^{5} + q^{4})s_{221}  + (q^{9} + q^{8} + q^{7} + q^{6})s_{2111}  + q^{10}s_{11111},\\
&W_{41}=s_{41} + (q^{2} + q)s_{32}+ (q^{3} + q^{2} + q)s_{311} + (q^{4} + q^{3} + q^{2})s_{221} + (q^{5} + q^{4} + q^{3})s_{2111}   + q^{6}s_{11111},\\
&W_{32}=s_{32}+ qs_{311}+ (q^{2} + q)s_{221}+ (q^{3} + q^{2})s_{2111}   + q^{4}s_{11111},\\
&W_{311}=s_{311}+ qs_{221} +(q^{2} + q)s_{2111} +  q^{3}s_{11111} ,\\
&W_{221}=s_{221}+ qs_{2111} + q^{2}s_{11111} ,\\
&W_{2111}= s_{2111}+qs_{11111},\\
&W_{11111}=s_{11111}.
\end{aligned}$}

\noindent
 One may check that, whenever $\mu$ is larger than $\nu$ in dominance order then $W_\mu-q^{\eta(\mu')-\eta(\nu')}W_\nu$ is Schur positive.

%%%%%%%%%%%%%%%%%%
%              Section                          %
%%%%%%%%%%%%%%%%%% 
\section{\bleu{Cauchy kernel and scalar product}}
  The Cauchy-kernel formula for the $q$-Whittaker polynomials may be coined as
\begin{equation}\label{Whittaker_Cauchy}
    h_n\left[\frac{\bm{x}\cdot \bm{y}}{(1-q)}\right] = \sum_{\mu\vdash n} W_\mu(\bm{x})\widehat{W}_\mu(\bm{y})
 \end{equation}
If one considers the scalar product $\langle -,-\rangle_q$ characterized on the power sum basis by: 
    $$\langle p_\lambda,p_\mu/z_\mu\rangle_q := \delta_{\lambda,\mu} \prod_{k\in\mu} (1-q^k),$$
then, the Cauchy-kernel formula is equivalent to 
\begin{equation}\label{dual_basis_Whittaker}
  \langle W_\lambda,\widehat{W}_\mu\rangle_q=\delta_{\lambda,\mu}.
\end{equation}
This $q$-scalar product is linked to the Hall scalar product via the relation
\begin{equation}\label{link}
   {\langle f,g\rangle = \langle f,g^\star\rangle_q},\qquad{\rm where}\qquad g^\star(\bm{z}):=g\Big[\frac{\bm{z}}{1-q}\Big].
 \end{equation}
It follows that the $q$-adjoint of skewing\footnote{Recall that $g^\perp$ is characterized by $\langle g^\perp f_1,f_2\rangle =\langle f_1,g\cdot f_2\rangle$, for the usual scalar product.} by $g$ is multiplication by $g^\star$:
  \begin{equation}\label{adjointness}
   \langle g^\perp f_1,f_2\rangle_q =\langle f_1,g^\star\!\cdot\!  f_2\rangle_q ,\qquad \hbox{for all}\qquad f_1,f_2.
 \end{equation}
 Observe that the Cauchy-kernel may then be written as $h_n^\star[\bm{x}\bm{y}]$, and that says that $\widehat{W}_{n}(q;\bm{z})=h_n^\star(\bm{z})$.
 
 \subsection*{Aspects of RSK} Taking the (Hall) scalar product of the Cauchy-Kernel for the Whittaker with $h_1(\bm{y})^n$ (this corresponds to taking the coefficient of $y_1y_2\cdots y_n$), one obtains
		$$\Big\langle h_1(\bm{y})^n, h_n\Big[\frac{\bm{xy}}{1-q}\Big] \Big\rangle= \frac{h_1^n(\bm{x})}{(1-q)^n}$$
and on the right-hand side:
$$\langle h_1^n(\bm{y}),  \sum_{\lambda\vdash n}W_\lambda(\bm{x})\widehat{W}_\lambda(\bm{y})\rangle=
         \sum_{\lambda\vdash n}W_\lambda(\bm{x})\, \langle h_1^n(\bm{y}),\widehat{W}_\lambda(\bm{y})\rangle  $$
with the convention that $\widehat{W}_\lambda$ stands for $W_\lambda$ divided by the product that we have had all along. Now, one recalls that, for any symmetric function $g$ of degree $n$, the scalar product $\langle g, h_1^n\rangle$ is equal to the iterated skewing $(h_1^\perp)^n g$.   For example, using in~\pref{Comb_Formula} the fact that $\langle s_\mu, h_1^n\rangle$ is the number $f^\lambda$ of standard Young tableaux, we get
	$$ \left(\frac{h_1}{1-q}\right)^3 = \frac{1}{(1-q)^3} W_3 + \frac{q+2}{(1-q)^2} W_{21} +\frac{1}{(1-q)} W_{111} $$
In general, we have (see~\pref{formula_g_mu})
\begin{equation}
      \left(\frac{h_1}{1-q}\right)^n=\sum_{\mu\vdash n} \frac{g_\mu(q)}{(1-q)^{\ell(\mu')}}\, W_\mu.
\end{equation}

%%%%%%%%%%%%%%%%%%
%              Section                          %
%%%%%%%%%%%%%%%%%%
\section{\bleu{\texorpdfstring{$W$}{W}-expansion of Macdonald polynomials}} 
The (modified) Macdonald polynomials have interesting expansion in the basis of the Whittaker polynomials (or Hall-Littlewood). To discuss this, let us consider the formula
\begin{equation}
	\mathcal{A}_\mu(q,t;\bm{z}):=\sum_{\nu} \tqbinom{n-\mu_1}{n-\nu_1}C_\mu^\nu(q,t)\ W_\nu(q;\bm{z}),
\end{equation}
where one sets $C_\mu^\nu(q,t)=0$ if $\nu\not\succeq \mu$, and 
\begin{equation}
   C_\mu^\nu(q,t):=\prod_{\myatop{c\in\mu\cap\nu}{\ell_\nu(c)\not=0}} (t^{\ell_\mu(c)}-q^{a_\mu(c)+1}),
\end{equation}
when $\nu\succeq \mu$, {\sl i.e.} $\nu$ is larger than $\mu$ in \define{dominance order}.
\begin{prop}
 In the special case when $\mu=(a,b)$ has at most two parts, with $n=a+b$, we have the equality
\begin{equation}\label{Hab_in_W}
	\MacH_{ab}(q,t;\bm{z})=\mathcal{A}_{ab}(q,t;\bm{z}).
\end{equation}
\end{prop}
Further formulas along these lines are:
\begin{align} 
	&\MacH_{k11}(q,t;\bm{z})=(t-qt^2)W_{(k+1,1)}+\mathcal{A}_{k11}(q,t;\bm{z});\\[4pt]
	&\MacH_{k21}(q,t;\bm{z})=
            (t-q^2)\,(t-q^{k-1})\,(t^2-q^k)\,W_{k21}(q;\bm{z})
          \nonumber   \\&\qquad\qquad\qquad\qquad
           +(1-q^{k-2})\,(t-q^{k-1})\,(t^2-q^k)\,W_{k3}(q;\bm{z})\nonumber \\
           &\qquad\qquad\qquad\qquad
           +(t-q^2)\,(t^2-q^k)\,W_{k+1,11}(q;\bm{z})\nonumber \\
           &\qquad\qquad\qquad\qquad
           +(t^2\,(q+t+1)+q^{2\,k-1}\,(q^2+q+1)\nonumber \\
           &\qquad\qquad\qquad\qquad\qquad\qquad-q^{k-1}\,(q+t)\,(q\,t+q+t))\,W_{k+1,2}(q;\bm{z})\nonumber \\
           &\qquad\qquad\qquad\qquad
           +(t\,(q+t+1)-q^k\,(q^2+q+1))\,W_{k+2,1}(q;\bm{z})
           +W_{k+3}(q;\bm{z}).
\end{align}
In general, the coefficient of $W_\mu$ in $\MacH_\mu$ is $C_\mu^\mu(q,t)$, and we thus have
\begin{prop} For all partitions $\mu$,
\begin{equation}\label{W_mellit}
\frac{1}{C_\mu^\mu(q,t)} \MacH_\mu(q,t;\bm{z}) = W_\mu(q;\bm{z})+\sum_{\nu\succ \mu} a_\lambda(q,t) W_\nu(q;\bm{z}).
\end{equation}
   \end{prop}
There is a similar identity for the Hall-Littlewood\footnote{This may readily be obtained from~\pref{W_mellit} via the ``Flip''.} which is exploited by Mellit in his paper~\cite[Cor. 5.12]{mellit}.

%%%%%%%%%%%%%%%%%%
%              Section                          %
%%%%%%%%%%%%%%%%%% 
\section{\bleu{Pieri rules}}

The  Pieri rules, have the following explicit expressions for the Whittaker. The ``adjoint''-Pieri rule corresponds to skewing by $h_k$, and we have the formula
\begin{equation}\label{adjoint_pieri_formula}
	h_k^\perp W_\lambda=\sum_{\mu\rightarrow_k \lambda} c_{\mu\lambda}(q) \,W_{\mu},
\end{equation}
with $c_{\mu\lambda}(q)$ given by formula~\ref{c_mu_nu}, and
where $\mu\rightarrow_k \lambda$ means that $\mu$ is obtained from $\lambda$ by removing a length $k$ horizontal strip. It follows that we have the recurrence (see~\cite{borodin_petrov})
\begin{equation}
   W_\lambda(q;\bm{z}+y)= \sum_{k\geq 0} \sum_{\mu\rightarrow_k \lambda} c_{\mu\lambda}(q)\, y^k\,W_\mu(q;\bm{z}),
\end{equation}
which allows the explicit calculation of $W_\mu$ on a given set of $m$-variables. For example,
\begin{align*}
   W_{42}(q;\bm{z}+y)&=W_{42}(q;\bm{z}) +y\left( (q  +1)W_{32}(q;\bm{z})  + (q +1)W_{41}(q;\bm{z}) \right) \\ 
     &\qquad  + y^2\left(W_{22}(q;\bm{z}) +(q+ 1)^2W_{31} (q;\bm{z}) + W_{4}(q;\bm{z}) \right) \\
     &\qquad  +y^3\left((q  + 1)W_{21}(q;\bm{z})   + (q + 1)W_{3} (q;\bm{z}) \right)+y^{4}\,W_{2} (q;\bm{z}) .
\end{align*}
Next, in view of~\pref{adjointness}, it is natural to consider that the $q$-Pieri rule is multiplication by the symmetric function $h_k^\star(\bm{z}) =h_k[\bm{z}/(1-q)]=\widehat{W}_{k}(q;\bm{z})$, in view of~\pref{formule_W_n}. Thus, using~\pref{dual_basis_Whittaker} we see that
   \begin{equation}
     c_{\mu\lambda}(q)= \langle h_k^\perp W_\lambda,\widehat{W}_\mu\rangle_q
                                  = \langle  W_\lambda,\widehat{W}_k\widehat{W}_\mu\rangle_q,
   \end{equation}
from which it follows that
  \begin{equation}\label{W_dual_Pieri}
     \widehat{W}_k\widehat{W}_\mu = \sum_\lambda c_{\mu\lambda}(q) \widehat{W}_\lambda.
    \end{equation}
Equivalently, since $\widehat{W}_\mu=W_\mu/v_\mu(q)$, 
 \begin{equation}\label{W_Pieri}
	     \widehat{W}_{k}\,W_\mu =\sum_{\mu\rightarrow_k \lambda}  d_{\mu\lambda}(q)  W_\lambda,\qquad {\rm where}\qquad d_{\mu\lambda}(q)=\frac{v_\mu(q)}{v_\lambda(q)}\,c_{\mu\lambda}(q),
	\end{equation}
so that we get the formula
 \begin{equation}\label{d_mu_nu}
   d_{\mu\lambda}(q)=\prod_{i\geq1} \qbinom{\sigma_\mu(i-1)}{\lambda_i-\mu_i}.
 \end{equation}
For instance, 
        $$\widehat{W}_{1}\,W_{52}=  (q + 1)W_{521} + (q^{2} + q + 1)W_{53} + \frac{1}{1-q}\,W_{62}$$

%%%%%%%%%%%%%%%%%% 
\subsection*{Combinatorial description of Pieri rules}
Simply put, the effect of $h_k^\perp$ on $W_\mu$ corresponds to removing from the diagram of $\mu$ $k$ cells among those that have leg-length equal to $0$.  In English notation (for those that do not speak French yet), these are the cells lying at the bottom (resp. top in French) of columns of the diagram. The weight of such a pierced diagram is $q^d W_\nu$, where $d$ is the number of left slide needed to turn this pierced diagram into a partition $\nu$, and one adds up all these weights over all possible configurations.
The analogous up-going Pieri rule, which corresponds to multiplication by $h_k[X/(1-q)]$, is similarly described.  Ties to representation theory are described in an upcoming section.

%%%%%%%%%%%%%%%%%%
%              Section                          %
%%%%%%%%%%%%%%%%%%
\section{\bleu{General up and down operators}} 
For any cell $c\in\N\times \N$, it is natural\footnote{See Proposition~\ref{Pieri_down_module}.} to consider \define{cell-indexed down operators}, $D_c$, (resp. \define{cell-indexed up operators}, $U_c$) recursively defined below on Whittaker polynomials. 
\begin{definition}[Down operators] The linear operator $D_c$, is such that
 for any partition $\lambda$, and any cell $c=(i,j)$,
 \begin{equation}\label{four_term_rec_down}
   D_c W_{\lambda}=\begin{cases}
     D_{c+(0,1)} W_{\lambda}+D_{c+(1,0)}W_{\lambda}-D_{c+(1,1)}W_{\lambda}, 
                 & \text{\rm if}\ \ell(c)>0\ {\rm with}\ c\in \lambda, \\[4pt]
     (1+\ldots +q^{a(c)})W_{\lambda\setminus \{c_{j+1}\}},&  \text{\rm if}\ \ell(c)=0\ {\rm with}\ c\in \lambda,\\[4pt]
        0,& \text{\rm otherwise},
\end{cases}
\end{equation}
where $c_k$ denotes the $k^{\rm th}$-row  inner corner {\rm (}if any{\rm )} of $\lambda$. 
 \end{definition}
We may check that the ``usual'' Whittaker $q$-Down operator corresponds to the case $c=(0,0)$, and that we get the general formula of Proposition~\ref{Dc_Pieri} below. 

\begin{prop}\label{Dc_Pieri}
For any partition $\lambda$, and any cell $c=(i,j)$, we have
\begin{equation}
 D_{c} W_\lambda= \sum_{c\leftmost_\mu c'}(1+q+\ldots+q^{a(c')})W_{\mu}.
 \end{equation}
In particular, it follows that
 \begin{equation}
 D_{(0,0)} W_\lambda=h_1^\perp W_\lambda.
 \end{equation}
In other words $D=D_{(0,0)}$.
\end{prop}
Next, let us set $\N^*:=\N\cup \{\infty\}$, with the addition rule $\infty\pm k=\infty$. We formally consider that
   $$1+q+\ldots+q^{\infty}=\frac{1}{1-q}.$$
Then in an approach similar to that above, we may define $U_c$ as follows.
\begin{definition}[Up operators] 
 For any partition $\mu$, and any $c=(i,j)$ in $\N^*\times\N^*$, we set
\begin{equation}
   U_c W_{\mu} :=\begin{cases}
     \sum_{c'\rightmost_\mu c} U_{c'}W_{\mu} & \text{\rm if}\ \ell(c)>0\ {\rm and}\ c\notin \mu,\\
     (1+q+\ldots +q^{a(c)})\,U_{c_k}W_{\mu} &\text{\rm if}\ \ell(c)=0\ {\rm and}\ c\notin \mu,\\ 
      0 & \text{\rm otherwise},
\end{cases}
 \end{equation}
 where $c_k$ denotes the $k^{\rm th}$-row outer corner {\rm (}if any{\rm )} of $\mu$. 
\end{definition}
The ``usual'' $q$-Whittaker Up operator corresponds to the case $c=(\infty,\infty)$, that is $U=U_{(\infty,\infty)} $: 
 \begin{equation}
 U_{(\infty,\infty)} W_\lambda=\widehat{W}_1\, W_\lambda.
 \end{equation}
From all this, we directly deduce a simple combinatorial understanding of the $q$-commutation rule 
   $$D\,U-U\,D=\frac{1}{1-q}\,\Id.$$
 
%%%%%%%%%%%%%%%%%%
%              Section                          %
%%%%%%%%%%%%%%%%%%
\section{\bleu{\texorpdfstring{$W$}{W}-positivity}} 
There are several interesting symmetric functions that are $W$-positive ({\sl i.e.} with coefficients in $\N[q,t]$). We present below three types of families of such symmetric functions.

%%%%%%%%%%%%%%%%%%
\subsection*{\underline{The elliptic Hall algebra paradigm}}
We consider here, the $W$-expansion of specializations of symmetric functions $Q_{mn}(q,t;\bm{z})$ that arise via the elliptic Hall algebra paradigm (see appendix).
Among the known identities between these, we recall that
\begin{equation}\label{skew_identity}
   e_1^\perp Q_{(n-1,n)}(q,t;\bm{z}) = Q_{(n-1,n-1)}(q,t;\bm{z}).
\end{equation}
 We also recall that, for the special cases $m=n-1$ and $m=n+1$, we have 
   $$Q_{(n-1,n)}(q,t;\bm{z})=\nabla((-qt)^{1-n} h_n)(q,t;\bm{z}),\qquad {\rm and}\qquad 
       Q_{(n+1,n)}(q,t;\bm{z})=\nabla(e_n)(q,t;\bm{z}).$$
This last symmetric function is the bigraded Frobenius characteristic of the space of diagonal $\S_n$-harmonics, aka the diagonal $\S_n$-coinvariant module. Hence, its specialization at $t=0$ is 
   \begin{equation}
       Q_{(n+1,n)}(q,0;\bm{z})=W_{n}(q;\bm{z}), 
    \end{equation}
 since this is the graded Frobenius characteristic of the classical  $\S_n$-coinvariant module (a.k.a. the cohomology ring of the full flag manifold). In fact, for all $m$, the specializations $Q_{mn}(q,0;\bm{z})$ are $W$-positive. When $\gcd(m,n)=1$, we have 
        $$Q_{mn}(q,0;\bm{z})=q^\beta W_{(m^k,r)}(q;\bm{z}),$$ 
 writing $(m^k,r)$ for the partition having $k$ parts of size $m$ and one part of size $r$, with $n=m\,k+r$ corresponding to Euclidean division. The parameter $\beta$ is simple to calculate, and often equal to $0$.  Some examples are:
{\tiny\begin{align*}
& Q_{19}(q,0;\bm{z}) =W_{1^9}(q;\bm{z}),  && Q_{29}(q,0;\bm{z}) =W_{2^41}(q;\bm{z}), && Q_{39}(q,0;\bm{z}) =(q + 1)W_{3321}(q;\bm{z}) + W_{3^3}(q;\bm{z}),\\
& Q_{49}(q,0;\bm{z}) =W_{4^21}(q;\bm{z}) ,&& Q_{59}(q,0;\bm{z}) =W_{54}(q;\bm{z}), && Q_{69}(q,0;\bm{z}) =(q^{3} + q^{2} + q)W_{63}(q;\bm{z}) ,\\
& Q_{79}(q,0;\bm{z})=q^{2}W_{72}(q;\bm{z}) && Q_{89}(q,0;\bm{z}) =W_{81}(q;\bm{z}).
\end{align*}}
In general, when $m=n$, we have
\begin{equation}
    Q_{nn}(q,0;\bm{z}) = W_{n}(q;\bm{z})+[n-1]_q\,W_{(n-1,1)}(q;\bm{z}).
\end{equation}
Also $W$-positive are the specialization at  $t=1/q$ of the $Q_{mn}$.  Indeed, using~\pref{formula_hook2}, we calculate that
 \begin{align}
q^{\alpha}Q_{mn}(q,1/q;\bm{z})&=\frac{[d]_q}{[m]_q}e_n[[m]_q\,\bm{z}]\nonumber \\
&=\frac{[d]_q}{[m]_q}e_n\!\left[{\textstyle \frac{(1-q^m)}{1-q}}\,\bm{z}\right]\nonumber\\
&=\frac{[d]_q}{[m]_q}\sum_{\mu\vdash{n}} \omega \widehat{W}_\mu[1-q^m]\,W_\mu(q;\bm{z})\nonumber\\
&= \sum_{\mu\vdash n}\Big(q^{\eta(\mu)} \frac{[d]_q}{[m]_q}
  \prod_{{\ell(i,j)=0}} \frac{1-q^{m-i}}{1-q^{a(i,j)+1}}\Big)W_\mu(q;\bm{z}),\label{Formula_Qnn}
\end{align}
where $d$ stands for ${\rm gcd}(m,n)$, and $\alpha:=(mn-m-n+d)/2$.  Some of the values of $q^{\alpha}Q_{mn}(q,1/q;\bm{z})$ are given in the appendix.

%%%%%%%%%%%%%%%%%
\subsection*{\underline{Specialization at \texorpdfstring{$(t=1/q)$}{t} of Delta operators}}
We may generalize the calculation in~\pref{Formula_Qnn}, to get the $W$-positivity of the specialization at $t=1/q$ of symmetric functions involved in the Delta-Conjecture. Let us start with the following identity (see \cite[Thm 5.1]{DeltaConjecture}), which holds for any degree $k$ homogenous symmetric function $f$: 
\begin{equation}
    q^{k\,(n-1)}\Delta_{f}(e_n)(q,1/q;\bm{z})=\frac{f[n]_q}{[k+1]_q}e_n([k+1]_q\bm{z})),
 \end{equation}
 where $f[n]_q$ stands for the plethystic evaluation $f[1+q+\ldots +q^{n-1}]$. For the definition of the Macdonald eigenoperator $\Delta_f$, see the appendix.
The Delta-Conjecture has to do with an explicit combinatorial formula for $\Delta_{e_k}(e_n)$, with $1\leq k\leq n$. If $f$ is any Schur positive symmetric function, then $f[n]_q$ lies in $\N[q]$. 
 
 For any symmetric function $\varphi(q;\bm{z})$ with coefficients in the fraction field $\Rational(q)$, let us write $\overline{\Delta}_f(\varphi)$ for the specialization\footnote{The eigenfunctions of the operator $\overline{\Delta}_f$ are the functions $s_\mu[\bm{z}/(1-q)]$, for $\mu$ any partition. } $(\Delta_{f}\,\varphi(q;\bm{z}))\big|_{t=1/q}$.
 We can imitate the calculation in~\pref{Formula_Qnn} to show that:
 
 \begin{prop}
  For all partitions $\mu$ of $k$, we have $W$-positivity of $q^{k (n-1)}\overline{\Delta}_{s_\mu}(e_n)$.
 \end{prop}

%%%%%%%%%%%%%%%%% 
\subsection*{\underline{Specialization at \texorpdfstring{$(t=0)$}{t} of Delta operators}}
Just as above, for any symmetric function $\varphi(q;\bm{z})$ with coefficients in the fraction field $\Rational(q)$, let us write $\Delta_{f}^\zero(\varphi)$ for the specialization $(\Delta'_{f}\,\varphi(q;\bm{z}))\big|_{t=0}$.
The eigenfunction of these operators are the Hall-Littlewood symmetric polynomials $\mathcal{H}_\mu(q;\bm{z})$,
with eigenvalue $f(q+\ldots+q^{\mu_1})$. Observe that we have the special case $\mathcal{H}_{n}=W_{n}$, hence $W_{n}$ is one of the eigenfunctions of $\Delta_{f}^\zero$. We also write  $\nabla^\zero$ for the specialization at $t=0$ of the operator $\nabla$. For all partitions $\mu\not=(n)$ of $n$, the operator $\nabla^\zero$ has eigenfunction $\mathcal{H}_\mu(q;\bm{z})$ with eigenvalue $0$; and for $\mathcal{H}_n(q;\bm{z})$ the eigenvalue is $q^{\binom{n}{2}}$. 

For any Schur positive $f$,  the symmetric function $\Delta^\zero_{s_\nu}(e_n)$ is $W$-positive. Indeed, Theorem 6.14 of~\cite{HRS1} may be formulated as:
\begin{align}\label{positivity_at_t0}
\Delta^\zero_{e_{k-1}}(e_n)
&=q^{-\binom{k}{2}}\sum_{\mu\vdash_k n}{q^{\eta(\mu')} \qbinom{k}{\bm{d}(\mu)} W_{\mu'}(q;\bm{z})},\\
&=q^{-\binom{k}{2}}\sum_{\mu\vdash n,\ \mu_1=k}{q^{\eta(\mu)} \qbinom{\mu_1}{\sigma(\mu)} W_{\mu}(q;\bm{z})}
\end{align}
 where we write $\mu\vdash_k n$ if $\mu$ is a length $k$ partition of $n$, and $\bm{d}(\mu)=(d_1,\ldots,d_n)$ when $\mu=1^{d_1}\cdots n^{d_n}$. This is to say that $d_i$ is the multiplicity of the part $i$ in $\mu$.
More generally, it is shown\footnote{However, our statements are expressed in terms of the symmetric functions $W_{\mu'}=\omega Q_\mu'$.} in~\cite[Thm 1.2]{HRS} that, for all partitions $\nu$ of $d$, we have $W$-positivity of $\Delta^\zero_{s_\mu}(e_n)$, with the explicit expansion
\begin{equation}\label{Delta_zero_Schur}
 \Delta^\zero_{s_\nu}(e_n)= \sum_{k=\ell(\nu)+1}^{|\nu|+1} P_{\nu,k-1}(q)
\sum_{\mu\vdash_k n}{q^{\eta(\mu')} \qbinom{k}{\bm{d}(\mu)} W_{\mu'}(q;\bm{z})}
\end{equation}
where
\begin{equation}
   P_{\nu,k}(q)= q^{d-k\,(k+1)} \sum_{\mu\vdash_k n} q^{\eta(\mu')} \qbinom{k}{\bm{d}(\mu)} K_{\nu,\mu}(q).  \qquad ({\rm where}\ \nu\vdash d)
 \end{equation}
Observe, in view of~\pref{W_Kostka} and~\pref{positivity_at_t0}, that $P_{\nu,k}(q)$ may readily be obtained (up to a power of $q$) by taking the scalar product  of $\Delta^\zero_{e_{k}}e_j$ with $s_{\nu'}(\bm{z})$. More precisely,
\begin{equation}
   P_{\nu,k}(q)=q^{d-\binom{k+1}{2}}\langle \Delta^\zero_{e_{k}} e_d,s_{\nu'}\rangle.
\end{equation}
Hence, for any degree $d$ homogeneous Schur positive symmetric function $f$, formula~\pref{Delta_zero_Schur} implies that we have $W$-positivity of
\begin{equation}
\Delta^\zero_{f}(e_n)= \sum_{k=0}^{d} q^{d-k} \langle \Delta^\zero_{e_{k-1}} e_d,\omega f\rangle\,\Delta^\zero_{e_{k}}(e_n),\qquad ({\rm where}\ d=\deg(f))
\end{equation}
In fact, this equality follows from the general operator identity
\begin{equation}
\Delta^\zero_f = \sum_{k=0}^{d} q^{d-k} \langle \Delta^\zero_{e_{k-1}}e_d,\omega f\rangle\,\Delta^\zero_{e_{k}}.
\end{equation}
Regarding this, it is interesting to recall (see~\cite[Lemma 6.1]{DeltaConjecture}) that for all $k$ and any homogeneous degree $d$ symmetric function $f$, one has
  $$\langle \Delta_{e_{k-1}}e_d,\omega f\rangle = \langle \Delta_f\, e_k,h_k\rangle.$$
Since $W_{n}$ is an eigenfunction of $\Delta^\zero_f$ (with eigenvalue $f(q+\ldots+q^{n-1})$), it follows that $\Delta^\zero_{e_k} W_{n}$ is trivially $W$-positive for all $k$. When $\Delta^\zero_{e_{k}}W_\mu$ is $W$-positive for all $k$, the above operator formula implies that $\Delta^\zero_{f}W_\mu$ is also $W$-positive for any Schur positive $f$.  Recall that $W_{1^n}=e_n$, so that formula~\pref{positivity_at_t0} has a bearing here.
One observes that
\begin{prop}
In the $W$-basis expansion of $\Delta^\zero_{e_k} W_\mu$, the coefficient of $W_\nu$ vanishes when $\nu$ is strictly dominated by $\mu$, or when $\nu_1<k+1$.
In particular, we have the explicit expansions
\begin{align}
  &\Delta^\zero_{e_{n-1}} W_\mu =q^{\eta(\mu')} W_{n}, &&\hbox{\rm for all } \mu\vdash n;\\
    &\Delta^\zero_{e_{n-2}} W_\mu =q^{\eta(\mu')+1-\mu_1}\tqbinom{\mu_1}{1}\,W_{n}+q^{\eta(\mu')}\tqbinom{n-\mu_1}{1}\,W_{(n-1,1)}, &&\hbox{\rm for all } \mu\vdash n;\\
&\Delta^\zero_{e_{k}} W_{ab} = \sum_{i=0}^b q^{\binom{k+1}{2}-i(k-b+1)}\tqbinom{b}{i}\tqbinom{a-1}{k-i}\,W_{(a+i,b-i)},
&&\hbox{\rm for all } a\geq b;\\
 &\Delta^\zero_{e_{k}} W_{(n-2,1,1)} =q^{\binom{k+1}{2}}\tqbinom{n-3}{k}\,W_{(n-2,1,1)}
 				+ q^{\binom{k}{2}}(q+1)\tqbinom{n-3}{k-1}\,W_{(n-1,1)}\nonumber \\
			&\qquad\qquad\qquad	+ q^{\binom{k-1}{2}}\tqbinom{n-3}{k-2}\,W_{n} ,
				&&\hbox{\rm for all }  n\geq 4.
\end{align}
Moreover, $\Delta^\zero_{e_k} W_\mu$ is $W$-positive (at least) for all $k$ and all two-part partitions and hook-shape partitions $\mu$.
\end{prop}
When $\mu\vdash n\leq 6$, the only case for which $\Delta^\zero_{e_k} W_\mu$ is not Schur positive (hence not $W$-positive either), is $\mu=222$, with $k=1$ or $2$. Indeed, we have
\begin{align*}
   \langle\Delta^\zero_{e_1}W_{222},s_{33}\rangle &=q^2-1,\\
    \langle\Delta^\zero_{e_2}W_{222},s_{33}\rangle &=q^6+q^5+2q^4+q^3-q;
 \end{align*}
and the corresponding $W$-expansions are
\begin{align*}
   \Delta^\zero_{e_1}W_{222} &= qW_{222} + \left(q^{3} + 2 q^{2} + q\right)W_{321} + \left(q^{2} - 1\right)W_{33},\\
   \Delta^\zero_{e_2}W_{222} &= \left(q^{4} + 2 q^{3} + q^{2}\right)W_{321} + \left(q^{3} - q\right)W_{33} + \left(q^{3} + q^{2}\right)W_{411} + \left(q^{5} + q^{4} + 2 q^{3}\right)W_{42}.
  \end{align*}

 %%%%%%%%%%%%%%% 
 \subsection*{\underline{Some \texorpdfstring{$(q,t)$}{qt}-cases of Delta operators}}
$W$-positivity also holds without restriction on $t$ in the following cases.
\begin{prop}
  For all $k$ and all partitions $\mu$ that have at most two parts, the symmetric function $\Delta'_{e_k}(W_\mu)$ is $W$-positive.
  We have the formulas
\begin{align}
  \Delta'_{e_{1}} W_{ab} &=\Delta^\zero_{e_1} W_{ab}+t\,[b]_q\,W_{ab},\\[4pt]
  \Delta'_{e_{2}} W_{ab} &=\Delta^\zero_{e_2} W_{ab}+(t\,q\,\tqbinom{a-1}{1}\tqbinom{b}{1}
                                             +t^2q\,\tqbinom{b}{2})\,W_{ab}
  					+t\,(q^b+q^{b-1})\tqbinom{b}{2}\,W_{(a+1,b-1)},\\
   \Delta'_{e_{3}} W_{ab} &= \Delta^\zero_{e_3} W_{ab}+(t\,q^3\tqbinom{a-1}{2}\tqbinom{b}{1}
   					+t^2\,q^2\tqbinom{a-1}{1}\tqbinom{b}{2}
   					+t^3\,q^3\tqbinom{b}{3}) W_{ab}\nonumber \\
				      &\qquad (t\,q^b\,[2]_q\tqbinom{b}{2}\tqbinom{a-1}{1}
				                        +t^2\,q^b\,[3]_q\,\tqbinom{b}{3})W_{(a+1,b-1)}\nonumber\\
				      &\qquad  +t\,q^{b+\binom{b-3}{1}}\,[3]_q\tqbinom{b}{3}\,W_{(a+2,b-2)}.
 \end{align}
 As well as (taking $k=n-1$)
\begin{align}
 \nabla W_{ab} &=q^{\textstyle \binom{a}{2}+\binom{b}{2}} \sum_{j=0}^b {\tqbinom{b}{j}}\,t^{b-j}\,W_{(a+j,b-j)}.
 \end{align}
\end{prop}

%%%%%%%%%%%%%%%
\subsection*{\underline{Weighted sums of LLT-polynomials}}
 A well-known combinatorial formula for $\nabla(e_n)$ may be expressed in terms of the LLT-polynomials $\mathbb{L}_\gamma(q;\bm{z})$ as
   $$\nabla(e_n) =\sum_\gamma t^{{\rm area}(\gamma)} \mathbb{L}_\gamma(q;\bm{z}),$$
 with $\gamma$ running over the set of $n$-Dyck paths. We consider the restriction of the right-hand side to Dyck paths to height at most $1$. 
 \begin{prop}
We have the positive $W$-expansion
 \begin{equation}
\mathcal{V}_n(q,t;\bm{z}):=  \sum_{{\rm height}(\gamma)\leq 1} t^{{\rm area}(\gamma)} \mathbb{L}_\gamma(q;\bm{z}) = 
  W_n+\sum_{j=0}^{\lfloor n/2\rfloor} \sum_{k=0}^n t^k \tqbinom{n-k}{j} \tqbinom{k-1}{j}\, W_{(n-j,j)}.
\end{equation}
\end{prop}
Recall that Dyck paths may be identified to partitions $\gamma$ contained in the staircase shape $(n-1,\ldots 2,1)$, and that the \define{area} and \define{height} of $\gamma$ are defined to be
$${\rm area}(\gamma):=\sum_{i} a_i, \qquad {\rm and}\qquad
   {\rm height}(\gamma):=\max_i \, a_i,$$
where the $a_i:=n-i-\gamma_i$ are the \define{row areas} of $\gamma$. We have the specialization
\begin{equation}
   \mathcal{V}_n(1,1;\bm{z})=\sum_{2a+b=n} \binom{n}{2a}\,e_{a^21^b}(\bm{z}),
\end{equation}
and the generating series expansion
\begin{equation}
  \frac{h_{1}(\bm{z}) -x t\,h_{2}(\bm{z})}{1 - x\,(t + 1)\,h_{1}(\bm{z}) + x^{2} t\,h_{2}(\bm{z})}=\sum_{n} \mathcal{V}_n(1,t;\bm{z})\,x^n.
\end{equation}
For $q=t=1$, the ``dimensions'' of these $\mathcal{V}_n$ is the number of forests of labeled trees, with at most one descent along any path from the root to a leaf (see~\cite{KK}).

%%%%%%%%%%%%%%%
\subsection*{\underline{\texorpdfstring{$\widehat{W}$}{W}-positivity}}
From the $\widehat{W}$-Pieri rule in~\pref{W_dual_Pieri}, we may derive the following.
\begin{prop}
For all partition $\mu$, we have\footnote{It is interesting to recall that $s_\mu\big[{\textstyle \frac{\bm{z}}{(1-q)}}\big]=
s_\mu\big[{\textstyle \frac{1}{(1-q)}}\big]\,\MacH_{\mu}(q,1/q;\bm{z})$.} 
   \begin{align}
       s_\mu\big[{\textstyle \frac{\bm{z}}{(1-q)}}\big] &=\det\left(\widehat{W}_{\mu_i+j-i}(q;\bm{z})\right)\\
       &=\widehat{W}_\mu(q;\bm{z})+\sum_{\lambda \succ \mu} \gamma_{\mu\lambda}(q)\,\widehat{W}_\lambda(q;\bm{z}),
   \end{align} 
      with the coefficients $\gamma_{\mu\lambda}(q)$ lying in $\N[q]$. When $\mu$ is a hook, then
      $$\gamma_{\mu\lambda}(q) =q^{w(\mu,\lambda)} \qbinom{\lambda_1-1}{\mu_1-1},$$
with $w(\mu,\lambda)=\eta(\lambda')-\binom{\lambda_1}{2}+\binom{\lambda_1-\mu_1+1}{2}$.
\end{prop}
This generalizes~\pref{formule_W_n},  as well as
\begin{align}
     e_n\big[{\textstyle \frac{\bm{z}}{(1-q)}}\big]   &=\left\langle e_n(\bm{x}), h_n\big[{\textstyle \frac{\bm{x}\bm{z}}{(1-q)}}\big]  \right\rangle\nonumber \\
    &= \sum_{\mu\vdash n} \big\langle e_n(\bm{x}),  W_\mu(q;\bm{x})\big\rangle\,\widehat{W}_\mu(q;\bm{z})\nonumber\\
    &=\sum_{\mu\vdash n} q^{\eta(\mu')}\widehat{W}_\mu(q;\bm{z}),
 \end{align}
 which is here derived from~\pref{Whittaker_Cauchy}. It is also worth recalling that $s_\mu\big[{\textstyle \frac{\bm{z}}{(1-q)}}\big]$ is an eigenfunction of $\overline{\nabla}$, with eigenvalue $q^{\eta(\mu')}$.

%%%%%%%%%%%%%%%%%%
%              Section                          %
%%%%%%%%%%%%%%%%%% 	
\section{\bleu{Graded \texorpdfstring{$\S_n$}{Sn}-modules for the Whittaker}}\label{sectionmodule}
Whittaker polynomials are closely related to the graded Frobenius characteristic of the cohomology ring of Springer fibers. One may actually get a representation theoretic description of $W_\mu$, in terms of a submodule of (generalized) Garsia-Haiman modules (see relevant appendix section) as follows. Using definition~\pref{M_mu_homog_component}, one sets
  \begin{equation}\label{Whittaker_Modules}
     \WN_\mu:=\bigoplus_{k} \M_\mu^{(k,\eta(\mu))},
   \end{equation}
 considering that the $\bm{y}$-variables are of degree $0$ (we say that they are \define{inert}). It follows from the $n!$-theorem that
 \begin{equation}\label{whittaker_partition}
 \WN_\mu(q;\bm{z})=W_\mu(q;\bm{z}).
 \end{equation}
 These modules are naturally ``anti-isomorphic'' (twisting by sign, and with a complement in degree), to the modules obtained by derivation closure of the classical realization of the Specht modules in the ring of polynomials in $\bm{x}$-variables. This is the $\bm{y}$-variable free part of $\M_\mu$, {\sl i.e.}
   \begin{equation}\label{HallLittlewood_Modules}
     \mathcal{S}_{\mu'}:=\bigoplus_{k} \M_\mu^{(k,0)}.
   \end{equation}
The anti-isomorphism map sends an element $f(\bm{x})\in \mathcal{S}_{\mu'}$ to $f(\partial \bm{x})\bm{V}_\mu(\bm{x},\bm{y})\in \WN_\mu$. In particular, this makes it clear that the degree zero component (do not forget that $\bm{y}$-variables are of degree $0$) of $\M_\mu$ is an irreducible representation corresponding to $s_\mu$. There is also an interesting inclusion of $\mathcal{S}_{\nu'}$ in $\mathcal{S}_{\mu'}$, when $\nu$ is dominated by $\mu$. This explains why one sees $q^{\eta(\mu')-\eta(\nu')}W_\nu$ lying inside $W_\mu$.
 
 %%%%%%%%%%%%%%%%
 \subsection*{Down-going Pieri}
 The following is deduced\footnote{By considering what it states about the highest $t$-degree components.} from the {\sl Four Term Recurrence Conjecture} of~\cite[Formula I.17]{poking_hole}, which is still open. However, in this instance we can prove all  statements. 
 Using the  generalized notion of Garsia-Haiman modules, associated to any determinant $\bm{V}_{\bm{d}}$ (see~\pref{defn_Delta_d}), we consider the case when $\bm{d}$ is a punctured Young diagram as follows. Let $\lambda$ be a partition of $n+1$, and consider one of the cells $c_i=(1+\lambda_{i+1},i-1)$ (in Cartesian coordinates for French, or matrix coordinates for English), on a row $i$ such that $\lambda_{i}>\lambda_{i+1}$.
\begin{figure}[ht]
\begin{center}
 \Yboxdim{13pt}
\setlength{\unitlength}{1em}
\begin{picture}(16,5)
\put(0,0){\gyoung(;;;;;;;;;;;;;,;;;;;;;;;;;;,;;;;;:{c_i}!\ygreen;;;;;,;;;;;,;;;;)}
\put(6.6,3.8){\put(0,0){$\flechegauche{1.5}$}
	\put(2.7,0){\small $a(c_i)$}
	\put(4.2,0){$\flechedroite{1.5}$}}
\put(-2,3.4){\put(1,0){$\flechehaut{3.2}$}
	\put(-1.8,0){$\flechebas{3.2}$}}
\put(-1.5,1.3){$i$}
\end{picture}
\end{center}
\qquad  \vskip-15pt
\caption{Punctured diagram $\lambda\setminus\{c_i\}$.}
\label{Fig_punctured} 
\end{figure}
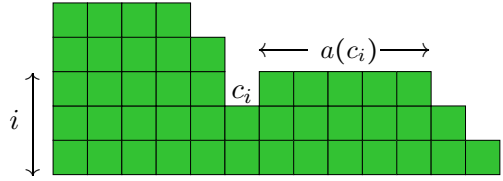
Inside the graded module $\M_{\lambda\setminus\{c_i\}}$, one then considers the top degree $\bm{y}$-component $\WN_{\lambda\setminus\{c_i\}}$ (as in~\pref{Whittaker_Modules}). 
\begin{prop}
With the degree of the variables $\bm{y}$ equal to zero, the graded Frobenius of $\WN_{\lambda\setminus\{c_i\}}$ is equal to
\begin{equation}\label{sub_corners}
   \WN_{\lambda\setminus\{c_i\}}(q;\bm{z})=(1+q+q^2+\ldots+q^{a(c_i)})\, W_\mu(q;\bm{z}),
 \end{equation}
where $\mu$ is the partition obtained by removing from $\lambda$ its corner in the $i^{\rm th}$-row. 
\end{prop}
Seeing this essentially reduces to the verification that the  operator  $\kappa_{\bm{x}}:=\sum_i {\partial x_i}$ on $\bm{V}_{\lambda\setminus\{c_i\}}$ has kernel of $\kappa_{\bm{x}}$ is $\WN_\mu$. Then one checks that the module 
$\WN_{\lambda\setminus\{c_i\}}$ decomposes into a direct sum of modules that are all isomorphic to $\WN_\mu$, suitably shifted in $\bm{x}$-degree. For example, we have:
\begin{center}
  \Yboxdim{7pt}
\setlength{\unitlength}{.5em}
\begin{picture}(35,7)
\put(-12,0){\gyoung(;;;;;;;;;;;;;,;;;;;;;;;;;;,;;;;;;:;!\ygreen;;;,;;;;;,;;;;)}
\put(7,2){\large$=(1+q+q^2+q^3)$}
\put(26,0){\gyoung(;;;;;;;;;;;;;,;;;;;;;;;;;;,;;;;;;;;;;!\ygreen,;;;;;,;;;;)}
\end{picture}
\end{center}
To make this more systematic, let us consider any cell $c$ in $\lambda$, and the corresponding module $\WN_{\lambda\setminus\{c\}}$.%, which we will here denote $\WN_ c$ (since $\lambda$ is fixed). 
We set $\WN_ {\lambda\setminus\{c\}}$ to be the trivial module $\{0\}$ if $c$ does not lie in $\lambda$.
As before, the associated graded Frobenius characteristic is $\WN_{\lambda\setminus\{c\}}(q;\bm{z})$. We have
\begin{prop}\label{Pieri_down_module}
 For any partition $\lambda$ and any cell  $c$ of $\lambda$, 
 \begin{equation}\label{four_term_rec_module}
   \WN_{\lambda\setminus\{c\}}(q;\bm{z})=D_c W_\lambda(q;\bm{z}).
\end{equation}
\end{prop}
Now, one may readily show that the module $\WN_{\lambda\setminus\{(0,0)\}}$ is isomorphic to the restriction to $\S_n$ of the $\S_{n+1}$-module $\WN_\lambda$. It thus follows that the down-going Pieri rule, {\i.e.} skewing by $h_1$,  may be interpreted as
    \begin{align}
        \WN_\lambda\Big\downarrow_{\S_n}^{\S_{n+1}} &\simeq \WN_{(0,0)}\\
             &\simeq\bigoplus_{\lambda_{i}>\lambda_{i+1}} \WN_{\lambda\setminus\{c_i\}}
      \end{align}
In terms of graded Frobenius characteristics, this corresponds to the case $k=1$ of~\pref{adjoint_pieri_formula}
   \begin{align*}\setlength{\fboxsep}{1pt}
   \WN_{\lambda\setminus\{(0,0)\}}(q;\bm{z})&=h_1^\perp W_\lambda(q;\bm{z})\\
    &=\WN_{\lambda\setminus\{c_i\}}(q;\bm{z})\\
    &=\sum_{i} [a(c_i)+1]_q W_{\setlength{\fboxsep}{1pt}\lambda^{\text{-}\framebox{$\scriptscriptstyle i$}}}(q;\bm{z}),  \end{align*}
where $\setlength{\fboxsep}{1pt}\lambda^{\text{-}\framebox{$\scriptscriptstyle i$}}$ denotes the partition obtained by removing the corner on row $i$ (the indices in the sum correspond to rows on which there is a corner). 
  In pictures, this is
\begin{center}
  \Yboxdim{7pt}
\setlength{\unitlength}{.5em}
\begin{picture}(35,14)     
       \put(0,8){\put(-15,0){\gyoung(:!\ygreen;;;;;;;;;,;;;;;;;;;;,;;;;;;;,;;;;)}
       \put(0,1.5){$=$}
       \put(5,0){\gyoung(;;;;;;;;;;,;;;;;;;:!\ygreen;;,;;;;;;;,;;;;)}
       \put(19.5,1.5){\large$+$}
       \put(23,0){\gyoung(;;;;;;;;;;,;;;;;;;;;;,;;;;:!\ygreen;;,;;;;)}
      \put(38,1.5){\large$+$}
       \put(41,0){\gyoung(;;;;;;;;;;,;;;;;;;;;;,;;;;;;;,:!\ygreen;;;)}}
       \put(0,1.5){{$=$}\ $[3]_q$}
       \put(6,0){\gyoung(;;;;;;;;;;,;;;;;;;;;:!\ygreen,;;;;;;;,;;;;)}
       \put(19.5,1.5){{\large$+$}\ $[3]_q$}
       \put(25,0){\gyoung(;;;;;;;;;;,;;;;;;;;;;,;;;;;;:!\ygreen,;;;;)}
      \put(38.5,1.5){{\large$+$}\ $[4]_q$}
       \put(44,0){\gyoung(;;;;;;;;;;,;;;;;;;;;;,;;;;;;;,;;;:!\ygreen)}
 \end{picture}
\end{center} 
On $\N\times \N$ we consider the component-wise partial order,  so that $(i,j)\preceq (k,l)$ if and only if $i\leq k$ and $j\leq l$.  Then, for $c$ any cell in $\mu$, we have
\begin{equation}\label{Generalized_Pieri_Down}
   \WN_{\lambda\setminus \{c\}} \simeq \bigoplus_{c\leftmost_\mu c'} \WN_{\lambda\setminus\{c'\}}.
\end{equation}

%%%%%%%%%%%%%%%%%
\subsection*{Up-going Pieri}Likewise, we have a representation theoretic description of the up-going Pieri rules (induction from $\S_{n-1}$ to $\S_n$), which corresponds to adding a box at the top of columns, including the $0$ columns sitting to the right of partitions (there are infinitely many of these). For $c\in\mu$, set  $\WN_{\mu\cup\{c\}}:=\{0\}$, then .
\begin{prop}
For any partition $\mu$, and $d$ a cell lying outside of $\mu$, we have a graded module isomorphism
\begin{equation}\label{Generalized_Pieri_Up}
   \WN_{\mu+ \{d\}} \simeq \bigoplus_{d'\rightmost_\mu d} \WN_{\mu+\{d'\}}.
 \end{equation}  
 Furthermore, the graded Frobenius characteristic of $\WN_{\mu+ \{d\}}$ is given by the formula
 \begin{equation}
    \WN_{\mu+ \{d\}}(q;\bm{z})=U_d\, W_\mu(q;\bm{z}).
\end{equation}
\end{prop}
In particular, adding to $\mu$ a cell at ``$(\infty,\infty)$'' corresponds to: either adding a cell anywhere on the first level to the right of the diagram, or adding a cell immediately above one of the corners of $\mu$. 
In pictures,  the above identity then takes the form
\begin{center}
  \Yboxdim{7pt}
\setlength{\unitlength}{.5em}
\begin{picture}(25,14)     
\put(3,8){  \put(-30,1.5){$W_{\mu+(\infty,\infty)}=$}  
       \put(-19,1.5){\Large$\frac{1}{1-q}$} \put(-15,0){\gyoung(;;;;;;;;!\yred;!\ygreen,;;;;;;;;,;;;;;,;;;)}
       \put(-2,1.5){\large$+$}
       \put(2,0){\gyoung(;;;;;;;;;,;;;;;;;;;,;;;;;;::!\yred;!\ygreen,;;;)}
       \put(15.7,1.5){\large$+$}
       \put(20,0){\gyoung(;;;;;;;;;,;;;;;;;;;,;;;;;;;,;;;:::!\yred;!\ygreen)}
      \put(34,1.5){\large$+$}
       \put(38,0){\gyoung(;;;;;;;;;,;;;;;;;;;,;;;;;;,;;;,::!\yred;!\ygreen)}}
       \put(-22.5,1.5){$=$}
       \put(-20,1.5){\Large$\frac{1}{1-q}$} \put(-15,0){\gyoung(;;;;;;;;;!\yred;!\ygreen,;;;;;;;;;,;;;;;;,;;;)}
       \put(-1.5,1.5){{\large$+$}\ $[3]_q$}
       \put(4.3,0){\gyoung(;;;;;;;;;,;;;;;;;;;,;;;;;;!\yred;!\ygreen,;;;)}
       \put(16.8,1.5){{\large$+$}\ $[4]_q$}
       \put(22.5,0){\gyoung(;;;;;;;;;,;;;;;;;;;,;;;;;;;,;;;!\yred;!\ygreen)}
      \put(35,1.5){{\large$+$}\ $[3]_q$}
       \put(41,0){\gyoung(;;;;;;;;;,;;;;;;;;;,;;;;;;,;;;,!\yred;!\ygreen)}
 \end{picture}
\end{center}     
Extension of this approach to general Pieri rules is closely tied to Theorem 2.3 of~\cite{multiple_left}, which may be stated with our current notations as
\begin{equation}
     \WN_{\bm{d}_{n,k}} (q;\bm{z})= \qbinom{n}{k}\, W_{k} (q;\bm{z}),
\end{equation}
for ${\bm{d}_{n,k}}:=\{(n-k,0),(n-k+1,0),\ldots,(n-1,0)\}$.

%%%%%%%%%%%%%%%%%%
%              Section                          %
%%%%%%%%%%%%%%%%%% 
\section{\bleu{Link to Science Fiction}}\label{section_SF}
The Whittaker expansion of the modified Macdonald polynomials $\MacH_\mu(q,t;\bm{z})$ is closely related to a new refinement of the ``Science-Fiction'' conjecture (which is described in~\cite{ScienceFiction}). This is related to the following relations among Whittaker polynomials and associated modules. To start our discussion of this, let us consider the polarization operator 
	$$\pi=y_1\partial x_1^r+y_2\partial x_2^r+\ldots +y_n\partial x_n^r,$$
as a $\S_n$-morphism between modules $\WN_{ab}$ indexed by two part partitions $\mu=(a,b)$, with $r=a-b+1$. Thus, we get a sequence of surjective morphism
  \begin{equation}\label{filtration}
  \xymatrix@-=0.6cm{
           \WN_{n}   \ar@{->>}[r]^-\pi 
	& \ar@{.}[r]
	&\ \ar@{->>}[r]^-\pi 
	& \WN_{ab}    \ar@{->>}[r]^-\pi  
	&\WN_{(a-1,b+1)}  \ar@{->>}[r]^-\pi  
	& \ar@{.}[r]
	&\ \ar@{->>}[r]^-\pi 
        & \WN_{k\ell} \ar@{->>}[r]^-\pi 
        &0,}
\end{equation}
with $k=\lceil n/2\rceil$ and $\ell=\lfloor n/2\rfloor$.
Respective kernels are denoted by $K_{ab}$, and we have graded $\S_n$-module isomorphisms
\begin{equation}
    \WN_{(a-1,b+1)} \simeq \WN_{ab}/K_{ab},\qquad \hbox{\rm (up to a degree shift by }a-b+1). 
 \end{equation}
More explicitly, the respective graded-Frobenius characteristics are related by the identity
\begin{align} 
   &(\WN_{ab}/K_{ab})(q;\bm{z}) = q^{a-b+1}W_{(a-1,b+1)}(q;\bm{z}),
  \end{align}
 and thus
 \begin{align}  
    & K_{ab}(q;\bm{z}) = W_{ab}(q;\bm{z})-q^{a-b+1}W_{(a-1,b+1)}(q;\bm{z}).
  \end{align}
We also observe that
  \begin{align}
       &K_{n}(q;\bm{z})=q^{\binom{n}{2}}\, W_{(n-1,1)}(1/q;\bm{z}). 
  \end{align}
We aim to exploit the above context to show and explain the identity~\pref{Hab_in_W}.

Let us next recall a construction of~\cite{ScienceFiction}. Given a partition $\lambda$ of $n+1$, let $\bm{c}$ be a subset of the partitions that may be obtained by removing a corner of $\lambda$, and let $0\leq j\leq m-1$, with $m$ standing for the cardinality of $\bm{c}$. As stated in~\cite[Thm 3.1]{ScienceFiction}) there are {\sl Schur positive}\footnote{The conjecture is about this positivity, and its interpretation in terms of submodules of the Garsia-Haiman modules.} symmetric functions $\Phi_{\bm{c}}^{(j)}$ (with a small change in the way indices are setup):
\begin{equation}
    \Phi_{\bm{c}}^{(j)}(q,t;\bm{z}):=\sum_{\nu\in C} P_\nu(\bm{c})\, (-T_\nu)^{j}\,  \MacH_\nu(q,t;\bm{z}),
    \qquad{\rm with}\qquad  P_\nu(\bm{c}):=\prod_{\alpha\in \bm{c}\setminus \{\nu\}}\frac{1}{1-T_\nu/T_\alpha};
\end{equation}
for which we have
   \begin{equation}
   \MacH_\nu(q,t;\bm{z})=\sum_{j=0}^{m-1} e_{j}[S_\nu({\bm{c}}) ]\ \Phi_{\bm{c}}^{(j)}(q,t;\bm{z}),
   \qquad{\rm with}\qquad S_\nu({\bm{c}}):=\sum_{\alpha\in \bm{c}\setminus \{\nu\}} \frac{1}{T_\alpha},
   \end{equation}
  whenever $\nu\in \bm{c}\subseteq \{\alpha\ |\ \alpha\rightarrow_1\mu\}$, for some partition $\mu$.
 Recall that we set $T_\mu:=\prod_{ab\in \mu} q^at^b$; and that this is the coefficient of $s_{1^n}$ in $\MacH_\mu$. When $\bm{c}$ is equal to 
 	$\{\mu\ | \mu\rightarrow_1\lambda\}$
is the full set $\{\mu\ | \mu\rightarrow_1\lambda\}$ of partitions that may be obtained from $\lambda$ by removing one of its corners, we also write $ \Phi_{\lambda}^{(j)}$ for $ \Phi_{\bm{c}}^{(j)}$. 

We aim to refine this, and to coin the Science Fiction conjecture in terms of the $W_\lambda$'s. As a first step, we consider partitions $\mu=(a,b)$ of $n$, so that $0\leq b\leq a$ with  $a+b=n$. For the rest of this section, let us set $k=\lceil n/2\rceil$ and $\ell=\lfloor n/2\rfloor$.
We now explain how all the polynomials $\MacH_{ab}$, decompose nicely as linear combinations (over $\N[q,t]$) of the $\N[q]$-symmetric function $\Cjn{i}{n}(q;\bm{z})$ that occur as coefficients (up to a $q$-binomial factor)  of $t^i$ in $\MacH_{k\ell}(q,t;\bm{z})$. More precisely, we claim the following.
\begin{prop} There are unique Schur positive symmetric functions such that
     \begin{equation}\label{Definition_C_jn}
       \MacH_{k\ell}(q,t;\bm{z}) = \sum_{j=0}^\ell \textstyle t^j\,\qbinom{\ell}{j}\,\Cjn{j}{n}(q;\bm{z}),
    \end{equation}
The $\Cjn{j}{n}$'s are graded Frobenius characteristics of modules that have dimension $n!/2^\ell$.
\end{prop}
 Observe that $C_\ell^n(q;\bm{z})=W_{k\ell}(q;\bm{z})$, and that we have the symmetry (implied by~\pref{Mac_flip})
      \begin{equation}\label{Symmetry_C_jn}
         q^d \Cjn{i}{n}(1/q;\bm{z}) = \Cjn{n-i}{n}(q;\bm{z}),
    \end{equation}
 where $d$ is the $q$-degree of $\Cjn{i}{n}$. 
For example, with $n=6$, we have the Schur expansions
\begin{align*}
    \Cjn{0}{6}&=q^{6}s_{222} + \left(q^{5} + q^{4}\right)s_{321} + q^{3}s_{33} + q^{3}s_{411} + \left(q^{4} + q^{3} + q^{2}\right)s_{42} + \left(q^{2} + q\right)s_{51} + s_{6},\\
    \Cjn{1}{6}&=q^{4}s_{2211} + q^{3}s_{3111} + \left(q^{3} + q^{2}\right)s_{321} + q^{2}s_{33} + \left(q^{2} + q\right)s_{411} + qs_{42} + s_{51},\\
    \Cjn{2}{6}&=q^{4}s_{21111} + q^{3}s_{2211} + q^{2}s_{222} + \left(q^{3} + q^{2}\right)s_{3111} + \left(q^{2} + q\right)s_{321} + qs_{411} + s_{42},\\
    \Cjn{3}{6}&=q^{6}s_{111111} + \left(q^{5} + q^{4}\right)s_{21111} + \left(q^{4} + q^{3} + q^{2}\right)s_{2211} + q^{3}s_{222} + q^{3}s_{3111} + \left(q^{2} + q\right)s_{321} + s_{33}.
\end{align*}
Thus, all the $\MacH_{(n-i,i)}$ can be expressed as a reunion of $2^\ell$ components, each of which  isomorphic to one of the $\Cjn{j}{n}$ (suitably degree shifted). For example, the matrices $(\gamma_{ij}^n)_{0\leq i,j\leq \ell}$ 
factorize nicely in LU-form as follows:
\begin{equation}
   \big(\gamma_{ij}^n\big)_{0\leq i,j\leq \ell} = \big(\tqbinom{i}{j} \prod_{a=0}^j (t-q^{n-i-a}) \big)_{ij}\, 
       \big(q^{(k-i)(j-i)} \qbinom{\ell-i}{j-i} \big)_{ij},
\end{equation}
with row and column indices starting at $0$. The first matrix (lower triangular) corresponds to~\pref{Hab_in_W}, and the second one (upper triangular) is the change of  basis between the $C_i^n$ and the $W_{(n-i,i)}$. In other words, 
\begin{prop}
We have the formula
\begin{equation} 
  W_{(n-i,i)}(q;\bm{z}) =  \sum_{j=0}^\ell q^{(k-i)(j-i)} \textstyle\qbinom{\ell-i}{j-i} C_j^n(q;\bm{z}),
 \end{equation}
where $n=k+\ell$, with $0\leq k-\ell\leq 1$. 
\end{prop}
One may readily invert this identity (since the transition matrix is upper uni-triangular), to
express the $C_i^n$ in terms of the $W_{(n-i,i)}$. Thus, we get
 \begin{equation}
    C_i^n(q;\bm{z}) =  \sum_{j=0}^\ell (-1)^{i+j}\tqbinom{\ell-i}{\ell-j}q^{\binom{k+1-i}{2}-\binom{k+1-j}{2}}\,W_{(n-j,j)}(q;\bm{z}).
 \end{equation}
For instance, with $n=6$ (hence $k=\ell=3$), we have 
 {\tiny \begin{align*}
\begin{pmatrix}
1 & q^{3}\,\qbinom{3}{1} & q^{6}\,\qbinom{3}{1} & q^{9} \\[3pt]
1 & t+q^{3} \qbinom{2}{1} & q^{6} + q^{2}t\,\qbinom{2}{1} & q^{4} t \\[3pt]
1 & q^{3} +t\, \qbinom{2}{1} & t^{2}+q^{2}t\, \qbinom{2}{1} & q t^{2} \\[3pt]
1 & t\,\qbinom{3}{1} & t^{2} \qbinom{3}{1} & t^{3}
\end{pmatrix}&=
\begin{pmatrix}
1 & 0 & 0 & 0\\[3pt]
1 & (t-q^{5}) & 0 & 0\\[3pt]
1 & \qbinom{2}{1}\,(t-q^{4}) & (t-q^{3})\,(t-q^{4}) & 0\\[3pt]
1 & \qbinom{3}{1}\,(t-q^{3}) & \qbinom{3}{1}\, (t-q^{2})\, (t-q^{3}) & (t-q)\,(t-q^{2})\,(t-q^{3})
\end{pmatrix}\times\\
&\hskip4cm \begin{pmatrix}
1& q^{3} \qbinom{3}{1}  &q^{6} \qbinom{3}{1} & q^{9}\\[3pt]
0&1& q^{2}\qbinom{2}{1}&q^{4}\\[3pt]
0&0&1&q\\
0&0&0&1
\end{pmatrix}
\end{align*}}
We also have the interesting formulas:
\begin{align}
 & \nabla\, \Cjn{0}{2n} (q;\bm{z}) = (-1)^n\,q^{n\,(3\,n-1)/2}\,t^n\,\Cjn{n}{2n}(q;\bm{z}), \qquad{\rm and}\\
  & \nabla\, \Cjn{0}{2n+1} (q;\bm{z}) = (-1)^n\,q^{3\,n\,(n+1)/2}\,t^n\,\Cjn{n}{2n+1}(q;\bm{z}).
\end{align}
Furthermore, the functions $\Phi_{(a+1,b)}^{(j)}$ have nice expansions in terms of the $\Cjn{j}{n}$ symmetric functions. For instance,
\begin{equation}
    \Phi_{(\ell+1,\ell)}^{(0)} = \sum_{i=0}^{\ell-1} \tqbinom{\ell-1}{i} \, t^i \, \Cjn{i}{2\ell}(q;\bm{z}).
 \end{equation}
We aim to refine this, and to coin the Science Fiction conjecture in terms of the $W_\lambda$'s. 

\section{Acknowledgments}
I would like to thank Florian Aigner, Benjamin Dequene, Gabriel Frieden, Steven Karp, Franco Saliola, and Hugh Thomas for many discussions, comments, and suggestions during the preparation of these notes. Part of section~\ref{section_SF} is the result of discussions with Jean-Christophe Novelli and Jean-Yves Thibon during my stay in the Institut Gaspard Monge in 2018-2019.

%%%%%%%%%%%%%%%%%%
%             Appendix                        %
%%%%%%%%%%%%%%%%%% 	
\setcounter{section}{1}
\renewcommand{\thesection}{\Alph{section}}

\section*{\bleu{Appendix}}
%%%%%%%%%%%%%%%%%% 
\subsection*{Symmetric functions and plethysm}\label{app_sym}
We mainly use Macdonald's notations (see~\cite{macdonald}).
In particular, we use the notation $\eta(\mu):=\sum_i(i-1)\mu_i$, for a partition $\mu=(\mu_i)$. 
For {cells} $(i,j)$ in a partition $\mu$, we often write $(i,j)\in\mu$.  We denote by $a_\mu(c)$ (resp. $\ell_\mu(c)$)  the arm length (resp. leg length) of the cell $c$ in a partition $\mu$. 

The usual {Hall scalar product} of $f$ and $g$ is denoted by $\scalar{ f}{g}$. The corresponding \define{classical Cauchy kernel} formula (which derives from RSK) may be expressed (using plethystic notation, see below) as
\begin{equation}\label{Classical_Cauchy}
    h_n[\bm{x}\bm{y}]=\sum_{\lambda\vdash n}s_\lambda(\bm{x})s_\lambda(\bm{y}).
 \end{equation}
 Here, one thinks of the set of variables $\bm{x}=(x_i)_i$ as a sum $\bm{x}=x_1+x_2+\ldots$ (likewise for $\bm{y}$). 
For a given symmetric function $g$, the linear operator $g^\perp$ is  the {adjoint} to the linear operation of multiplication by $g$. In formula, $\scalar{  g\cdot f_1}{f_2}=\scalar{  f_1}{g^\perp f_2}$,  for any symmetric functions $f_1$ and $f_2$. 
The classical {``adjoint'' Pieri rule}  (resp. its dual) (see~\cite{macdonald}), states that
  $$h^\perp_k s_\mu = \sum_{\nu\rightarrow_k \mu} s_\nu\qquad \hbox{(resp.}\ e^\perp_k s_\mu =\textstyle\sum_{\nu\rightarrow_k \mu} s_\nu)$$
with the indices of the sum running over partitions $\lambda$ that can be obtained from $\mu$ by removing $k$ cells, no two of which lying on the same column (resp. row).  The usual {Pieri rule} (resp. its dual) (see~\cite{macdonald}), says that
  $$h_k\, s_\mu = \sum_{\mu\rightarrow_k \lambda} s_\lambda\qquad \hbox{(resp.}\ e_k\, s_\mu =\textstyle \sum_{\mu\rightarrow_k \lambda} s_\lambda)$$
with the indices of the sum running over partitions $\lambda$ that can be obtained from $\mu$ by adding $k$ cells, no two of which lying on the same column (resp. row). Equivalently, the skew shape $\lambda/\mu$ is a horizontal band.

%%%%%%%%%%%%%%%%%% 
\subsection*{Plethysm}
For symmetric function $f$ and $g$, the \define{plethysm} $f\circ g=f[g]$ is a special case of $\lambda$-ring calculations $f[\bm{a}]$, in  which  symmetric function are considered as operators. The following evaluation rules entirely characterize these, assuming that $\alpha$ and $\beta$ are scalars, and that $\bm{a}$ and $\bm{b}$ lie in some suitable ring:
\begin{align}\itemsep=6pt
{\rm{(i)}}\quad &(\alpha f+\beta g)[\bm{a}]=\alpha\, f[\bm{a}]+\beta\, g[\bm{a}],\qquad &{\rm{(ii)}}&\quad (f\cdot g)[\bm{a}]=f[\bm{a}]\cdot g[\bm{a}],\nonumber\\
{\rm{(iii)}}\quad &p_k[\bm{a}\pm \bm{b}]=p_k[\bm{a}]\pm p_k[\bm{b}],                                   &{\rm{(iv)}}&\quad p_k[ \bm{a}\cdot \bm{b}]=p_k[\bm{a}]\cdot p_k[\bm{b}],\nonumber\\
{\rm{(v)}}\quad &p_k[ \bm{a}/\bm{b}]=p_k[\bm{a}]/p_k[\bm{b}],                                              &{\rm{(vi)}}&\quad p_k[p_j]=p_{kj},\label{def_plethysme}\\
{\rm{(vii)}}\quad &p_k[p_j\otimes p_\ell]=p_{kj}\otimes p_{k\ell},              &{\rm{(viii)}}&\quad p_k[\varepsilon] =(-1)^k,\nonumber\\
{\rm{(ix)}}\quad &p_k[x]=x^k,\ \hbox{if}\  x\  \hbox{is a ``variable''},                   
&{\rm{(x)}}&\quad p_k[c]=c,\ \hbox{if}\  c\  \hbox{is  a ``constant''}.\nonumber
\end{align}
\noindent In other words, for any ``atomic'' $v$, we must specify if $v$ is to be considered as a variable or a constant. The first two properties make it clear that any evaluation of the form $f[\bm{a}]$ may be reduced to instances of the form $p_k[\bm{a}]$. We also assume that property (iii) extends to denumerable sums. Property (viii) implies that $f[\varepsilon \bm{z}] = \omega\, f(\bm{z})$, for all symmetric function $f$. It may readily be shown that, for all symmetric function $f$,
\begin{equation}\label{e_skewing}
    \sum_{k\geq 0} u^k (e_b^\perp f)(q) = f[q-\varepsilon u].
 \end{equation}
See~\cite{bergeron,haglund} for more on plethysm.

%%%%%%%%%%%%%%%%%% 
\subsection*{Macdonald polynomials, and operators}\label{app_mac}
The (unmodified) Macdonald polynomials are denoted $P_\mu(q,t;\bm{z})$. They were originally as follows.
In his original paper of 1988 (see~\cite{macdonald_lotha}), Macdonald establishes  the existence and uniqueness of symmetric functions (polynomials) $P_\mu={P_\mu({\bf z};q,t)}$\index{Macdonald symmetric functions} such that
\begin{itemize}\itemsep=3pt\itemindent=20pt
  \item[(1)] $P_\mu=m_\mu+\displaystyle{\sum_{\lambda \prec \mu} c_{\mu \lambda}(q,t)}\, m_\lambda$,\quad
with coefficients $c_{\mu \lambda}(q,t)$ in $\Rational(q,t)$; and
  \item[(2)] $\langle P_\lambda,P_\mu \rangle_{q,t}=0$, whenever $\lambda \not= \mu$. 
 \end{itemize}
Recall that ${\lambda \prec \mu}$ stands for 
$\lambda$ being smaller than $\mu$ in the dominance order. 
In the above, the scalar product considered is characterized by:
\begin{equation}\label{scal_mac}
   {\langle p_\lambda,p_\mu \rangle}_{q,t}=\begin{cases}z_{\lambda}\,\textstyle  \prod_{i=1}^{\ell(\lambda)} \frac{1-q^{\lambda_i}}{1-t^{\lambda_i}} & \text{if}\ \lambda=\mu,\\[10pt]
                               0 & \text{otherwise},\end{cases}
\end{equation}
The \define{(pre)-modified Macdonald polynomials} 
    $$H_\mu(q,t;\bm{z}):= P_{\mu}\left[ \frac{{\bf z}}{1-t};q,t^{-1} \right]
                                  \prod_{c\, \in\, \mu} (q^{a(c)}-t^{\ell(c)+1})$$
where first briefly considered, and then they were replaced by the \define{modified Macdonald polynomials} $\MacH_\mu(q,t;\bm{z})$. 
The relation between the two is simply
   $$\MacH_\mu(q,t;\bm{z})= t^{\eta(\mu)}\,H_\mu(q,1/t;\bm{z}).$$
 It is a fact that $\eta(\mu)$ is the largest power of $t$ that occurs in both $H_\mu$ and $\MacH_\mu$. For $H_\mu$ (resp. $\MacH_\mu$), this is always in the coefficient of $s_{n}$ (resp. $s_{1^n}$).
 For example,
\begin{align*}
    H_{31}&=t\, s_{4} + (q^{2}t + qt + 1)\,s_{31}+ (q^{2}t + q )\,s_{22}  + (q^{3}t + q^{2}  + q )s_{211} + q^{3} s_{1111},\\  
    \MacH_{31}&=s_{4} + (q^{2} + q + t)\,s_{31}+ (q^{2} + q t)\,s_{22}  + (q^{3} + q^{2} t + q t)s_{211} + q^{3} t\,s_{1111}. 
 \end{align*}
 The coefficients of the expansion $H_\lambda(q,t;\bm{z})=\sum_\mu \widetilde{K}_{\lambda\mu}(q,t)\,s_\mu(\bm{z})$, are the $(q,t)$-Kostka polynomials. For example,
 $$(\widetilde{K}_{ \lambda\mu}(q,t))=\left(\begin{array}{rrrrr}
1 & q^{3} + q^{2} + q & q^{4} + q^{2} & q^{5} + q^{4} + q^{3} & q^{6} \\
t & q^{2} t + q t + 1 & q^{2} t + q & q^{3} t + q^{2} + q & q^{3} \\
t^{2} & q t^{2} + q t + t & q^{2} t^{2} + 1 & q^{2} t + q t + q & q^{2} \\
t^{3} & q t^{3} + t^{2} + t & q t^{2} + t & q t^{2} + q t + 1 & q \\
t^{6} & t^{5} + t^{4} + t^{3} & t^{4} + t^{2} & t^{3} + t^{2} + t & 1
\end{array}\right)$$
One sees here both the matrices $t^{\eta(\mu)}\,K_{ \lambda\mu}(1/t)$ and $K_{ \lambda'\mu'}(q)$ siting inside the matrix $(\widetilde{K}_{ \lambda\mu}(q,t))$, respectively setting $q=0$ and $t=0$.
There are analogous (modified) $(q,t)$-Kostka polynomials, such that $\MacH_\lambda(q,t;\bm{z})=\sum_\mu \widetilde{K}_{\lambda\mu}(q,t)\,s_\mu(\bm{z})$.
The set of {modified Macdonald polynomials}  
 $\{\MacH_\mu(q,t;\bm{z})\}_{\mu\vdash n}$
  forms a linear basis of the ring $\Lambda(q,t)$, of symmetric functions in the variables $\bm{z}=(z_i)_{i\in\N}$ over the field $\Rational(q,t)$. They are uniquely characterized by the equations 
 \begin{equation}
 \begin{array}{llll}
    \mathrm{(i)}\ {\displaystyle  \scalar{ s_\lambda(\bm{z})}{\MacH_\mu[q,t;(1-q)\,\bm{z}]}=0},\qquad{\rm if}\qquad  \lambda\not\succeq\mu,\\[6pt]
    \mathrm{(ii)}\ {\displaystyle   \scalar{ s_\lambda(\bm{z})}{ \MacH_\mu[q,t;(1-t)\,\bm{z}]}=0},\qquad{\rm if}\qquad   \lambda\not\succeq\mu',\ \mathrm{and}\\ [6pt]
    \mathrm{(iii)}\ {\displaystyle   \scalar{ s_n(\bm{z})}{ \MacH_\mu(q,t;\bm{z})}=1},
      \end{array}
 \end{equation} 
 involving plethystic notation.
See~\cite[Section 3.5]{haimanhilb} for more details on these dominance order triangularities. By definition, a Macdonald  \define{Eigenoperator} is a linear operator having the Macdonald polynomials  $\MacH_\mu$ as joint eigenfunctions. Hence, it is characterized by its eigenvalues on each of these. The Macdonald eigenoperator $\nabla$ has the  eigenvalues $T_\mu:=q^{\eta(\mu')}t^{\eta(\mu)}$ for $\MacH_\mu$, whereas the eigenvalues of the Macdonald operator $D_0$ are (see section~\ref{section_notations})
    $$1-(1-t)(1-q)B_\mu(q,t).$$
Some formulas pertaining to the $\MacH_\mu$ are
\begin{equation}
{ \MacH_\mu}(1,1;\bm{z})=h_1^n(\bm{z}),\qquad  \scalar{ \MacH_\mu}{s_{1^n}}  = T_\mu, \qquad {\rm and}\qquad \scalar{ \MacH_n}{s_{(a\,|\,b)}} = q^{\binom{b+1}{2}}\qbinom{n-1}{b}\label{Hn_binom}
\end{equation}
We also have the symmetries
   \begin{equation}\label{Mac_flip}
   \MacH_\mu(q,t;\bm{z}) =T_\mu\,\omega \MacH_\mu(1/q,1/t;\bm{z}), \qquad {\rm and}\qquad 
   \MacH_\mu(t,q;\bm{z}) =\MacH_{\mu'}(q,t;\bm{z}),
   \end{equation}   
which may readily be seen on the matrix $(\widetilde{K}_{\lambda\mu}(q,t))$.
We also recall that $\Delta_{f}$ (resp. $\Delta_{f}'$ ) is the Macdonald eigenoperator with eigenvalue $f[B_\mu(q,t)]$ (resp. $f[B_\mu(q,t)-1]$) for the eigenfunction $\widetilde{H}_\mu$. The case $f=e_n$ corresponds to $\nabla$ (as well as $m=n+1$ in \pref{Formula_Qnn}); and the case $f=e_1$ (up to a factor and shift by the Identity-operator) to the $\mathcal{O}_{10}$-operator of the elliptic Hall algebra (described further below). For more on operators on the (modified) Macdonald polynomials and their properties, see~\cite{identity}.

%%%%%%%%%%%%%%%%% 
\subsection*{Hall-Littlewood polynomials}
The (dual) Hall-Littlewood polynomials $\mathcal{H}_\mu(q;\bm{z})$ (also denoted $Q'_\mu(q;\bm{z})$) may be obtained by specializations from the Macdonald polynomials $\MacH_\mu(q,t;\bm{z})$. They occur as graded characters of singular cohomology ring of Springer fibers (see~\cite{de_concini_procesi,garsia_procesi}). They are related to the Schur functions by the following formula, involving the Kostka-Foulkes polynomials $K_{ \lambda\mu}(q)$  (which afford a combinatorial description in terms of the ``charge'' statistic on semi-standard Young tableaux). 
   \begin{equation}\label{Formula_HL}
       \mathcal{H}_\mu(q;\bm{z}) =\sum_\mu K_{ \lambda\mu}(q) s_\lambda(\bm{z}).
     \end{equation}
 For example, we have 
\begin{equation}\label{usual_kostka_matrix}
\Big(\langle\mathcal{H}_\mu,s_\lambda\rangle\Big)_{\mu,\lambda\vdash 4}=\left(\begin{array}{rrrrr}
1 & q^{3} + q^{2} + q & q^{4} + q^{2} & q^{5} + q^{4} + q^{3} & q^{6} \\
1 & q^{2} + q & q^{2} & q^{3} & 0 \\
1 & q & q^{2} & 0 & 0 \\
1 & q & 0 & 0 & 0 \\
1 & 0 & 0 & 0 & 0
\end{array}\right)
\end{equation}
%%%%%%%%%%%%%%%%% 
\subsection*{Cauchy Kernel and scalar product}
We set $$ w_\mu(q,t):=\prod_{c\in\mu}(q^{a(c)}-t^{\ell(c)+1})
    (t^{\ell(c)}-q^{a(c)+1}),$$
and then
    $$ \widehat{H}_\mu:=\frac{1}{w_\mu(q,t)} \MacH_\mu.$$
The Cauchy-kernel formula for the modified Macdonald states that
\begin{equation}
   e_n\!\left[\frac{\bm{x}\cdot\bm{y}}{(1-q)(1-t)}\right] = \sum_{\mu\vdash n} \MacH_\mu(\bm{x})\widehat{H}_\mu(\bm{y})
 \end{equation}
and the associated scalar product, denoted by $\langle -,-\rangle_\star$, is defined on the power sum basis as: 
   $$\langle p_\lambda,p_\mu/z_\mu\rangle_\star := (-1)^{n-\ell(\mu)} \delta_{\lambda,\mu} \prod_{k\in\mu} (1-q)^k(1-t^k).$$
 The Cauchy-kernel formula is equivalent to
 $$\langle \MacH_\lambda,\widehat{H}_\mu\rangle_\star=\delta_{\lambda,\mu}.$$
Observe that this scalar product is linked to the Hall scalar product via the relation
 \begin{equation}\label{starlink}
   {\langle f,g\rangle = \langle f,\omega\, g^\star\rangle_\star},\qquad{\rm where}\qquad g^\star(\bm{z}):=g\!\left[\frac{\bm{z}}{(1-q)(1-t)}\right].
 \end{equation}
Hence, the $\star$-adjoint of skewing by $g$ is multiplication by $\omega\, (g^\star)$.

%%%%%%%%%%%%%%%%%% 
\subsection*{Some explicit values are as follows.}\qquad

{\small
$\begin{aligned}
&\MacH_{2}=s_2+q\,s_{11},\\ 
&\MacH_{11}=s_2+t\,s_{11}; \\[6pt]
&\MacH_{3}=s_3+(q^2+q)\,s_{21}+q^3\,s_{111} ,\\ 
&\MacH_{21}=s_3+(q+t)\,s_{21}+q\,t\,s_{111} ,\\ 
&\MacH_{111}=s_3+(t^2+t)\,s_{21}+t^3\,s_{111};\\[6pt]
&\MacH_{4}=s_4+(q^3+q^2+q)\,s_{31}+(q^4+q^2)\,s_{2\,2}+(q^2+q+1)\,s_{211}\,q^3 +q^6\,s_{1111},  \\
&\MacH_{31}=s_4+(q^2+q+t)\,s_{31}+(q^2+q\,t)\,s_{2\,2}+q\,(q^2+q\,t+t)\,s_{211} +q^3\,t\,s_{1111},  \\
&\MacH_{22}=s_4+(q\,t+q+t)\,s_{31}+(q^2+t^2)\,s_{2\,2}+q\,t\,(q+t+1)\,s_{211} +q^2\,t^2\,s_{1111} ,\\
&\MacH_{211}=s_4+(t^2+q+t)\,s_{31}+(q\,t+t^2)\,s_{2\,2}+t\,(q\,t+t^2+q)\,s_{211} +q\,t^3\,s_{1111}, \\
&\MacH_{1111}=s_4+(t^3+t^2+t)\,s_{31}+(t^4+t^2)\,s_{2\,2}+t^3\,(t^2+t+1)\,s_{211}+t^6\,s_{1111}.
 \end{aligned}$}
 
%%%%%%%%%%%%%%%%
\subsection*{Elliptic Hall algebra} We start by defining operators $\mathcal{O}_{mn}$ on symmetric functions with coefficients in $\Rational(q,t)$. For more in this, see~\cite{elliptic}. The operator $\mathcal{O}_{01}$ is simply multiplication by the symmetric function $e_1(\bm{z})$, and $\mathcal{O}_{01}$ is the Macdonald eigenoperator $D_0$. For $m,n\geq 1$, we recursively define $\mathcal{O}_{mn}$ by the Lie bracket formula
\begin{displaymath}
         \mathcal{O}_{mn}:=\frac{1}{(1-t)(1-q)}\, [\mathcal{O}_{uv},\mathcal{O}_{k\ell}],
     \end{displaymath}
where $(k,\ell)$ is such that $(m,n)=(k,\ell)+(u,v)$, with $(k,\ell)$ and $(u,v)$ lying in $\N^2$,  $\ell-(k n/m)$ minimal, and 
\begin{displaymath}
   \det\begin{pmatrix} u & v\\ k& \ell\end{pmatrix}=d,
       \end{displaymath}
 with $d=\gcd(m,n)$. 
For $(m,n)=(ad,bd)$, we ask that $(k,\ell)$ be chosen to be the same as it would for $(a,b)$. 
 For example, we get the Lie bracket expressions
\begin{displaymath}
   \mathcal{O}_{43}=\frac{1}{(1-t)^6(1-q)^6} [[e_{{1}},D_{{0}}],[[e_{{1}},D_{{0}}],[[e_{{1}},D_{{0}}],D_{{0}}]]],
          \end{displaymath}
     or
\begin{displaymath}
   \mathcal{O}_{63}=\frac{1}{(1-t)^8(1-q)^8}[[e_{{1}},D_{{0}}],[[[e_{{1}},D_{{0}}],D_{{0}}],[[[e_{{1}},D_{{0}}],D_
{{0}}],D_{{0}}]]].
          \end{displaymath}
We then define the symmetric functions $Q_{mn}(q,t;\bm{z})$ to be the result of applying the operator $\mathcal{O}_{mn}$ to the symmetric function $s_0(\bm{z})=1$. Hence
\begin{equation}
    Q_{mn}(q,t;\bm{z}):=\mathcal{O}_{mn}(1).
 \end{equation}
For example, we have
\begin{align*}
Q_{13}&=s_{111},\\
Q_{23}&=\left(q + t\right)s_{111} + s_{21},\\
Q_{33}&=\left(q^{3} + q^{2} t + q t^{2} + t^{3} + q^{2} + 2 q t + t^{2} + q + t\right)s_{111} + \left(q^{2} + q t + t^{2} + 2 q + 2 t + 1\right)s_{21} + s_{3},\\
Q_{43}&=\left(q^{3} + q^{2} t + q t^{2} + t^{3} + q t\right)s_{111} + \left(q^{2} + q t + t^{2} + q + t\right)s_{21} + s_{3}.
\end{align*}
One may see that
\begin{equation}
   Q_{1,n}(q,t;\bm{z})=e_n(\bm{z}),\qquad {\rm and}\qquad \nabla\,Q_{mn}(q,t;\bm{z})=Q_{m+n,n}(q,t;\bm{z});
 \end{equation}
 as well as
 \begin{equation}
   Q_{nn}(q,t;\bm{z})= \nabla\big(\textstyle \sum_{k+\ell=n-1} (-qt)^{-k} s_{(k+1,1^\ell)}\big).
 \end{equation}

%%%%%%%%%%%%%%%%%% 
\subsection*{Bi-graded \texorpdfstring{$\S_n$}{Sn}-modules}
Each $\MacH_\mu$ polynomial affords an interpretation as the Frobenius transform\footnote{Which simply corresponds to naturally encode irreducible representations by Schur function.} of a (bi)-graded $\S_n$-modules (which is of linear dimension $n!$), known as the Garsia-Haiman module\footnote{See original papers~\cite{garsia_haiman_factorial,garsia_haiman_bigraded}.} (described below). These are (bi)-homogenous submodules $\M_\mu$ of the polynomial\footnote{In this section the variables $x_i$ and $y_i$ play a different role than in the rest of the text.} ring $\R=\Rational(\bm{x},\bm{y})$, in the variables $\bm{x}=x_1,\ldots x_n$ and $\bm{y}=y_1,\ldots y_n$, on which the symmetric group acts (diagonally)  by permutations of the variables: {\sl i.e.} $\sigma\cdot x_i=x_{\sigma(i)}$ and $\sigma\cdot y_i=y_{\sigma(i)}$. The grading is over $\N\times \N$, with value $(\deg_{\bm{x}}(f),\deg_{\bm{y}}(f))$ for $f\in\R$, and
\begin{equation}\label{GM_Modules}
    \M_\mu=\bigoplus_{k,j} \M_\mu^{(k,j)},
  \end{equation}
with $\M_\mu^{(k,j)}$ denoting the degree $(k,j)$  homogenous component of $\M$. Recall that the Frobenius transform of the graded character of such a module $\M_\mu$ may be defined as: 
   $$\M_\mu(q,t;\bm{z})=\sum_{k,j} q^kt^j \frac{1}{n!} \sum_{\sigma\in \S_n} \charac^{(k,j)}(\sigma) p_{\lambda(\sigma)}(\bm{z}),$$
 where $\charac^{(k,j)}$ stands for the character of $ \M_\mu^{(k,j)}$, and $\lambda(\sigma)$ is the partition giving the cycle decomposition of $\sigma$. The Schur expansion of $\M_\mu(q,t;\bm{z})$ gives the graded decomposition of $\M_\mu$ into irreducibles, which means that
   $$\M_\mu(q,t;\bm{z})=\sum_{k,j} q^kt^j \sum_{\nu\vdash n} c_{\nu}^{(k,j)} s_\nu(\bm{z}),\qquad {\rm iff}\qquad
     \M_\mu(q,t;\bm{z})\simeq \bigoplus_{k,j} \bigoplus_{\nu\vdash n} \mathcal{V}_\nu^{\oplus c_{\nu}^{(k,j)}},$$
where $\{\mathcal{V}_\nu\}_\nu$ constitute a complete family of representatives of the irreducible representations of $\S_n$, with  $c_{\nu}^{(k,j)}$
standing for the multiplicity of copies of $\mathcal{V}_\nu$ in $\M_\mu^{(k,j)}$.

%%%%%%%%%%%%%%%%%% 
\subsection*{Garsia-Haiman modules}
For each $n$-subset $\bm{d}$ of $\N\times \N$, called a diagram, one considers the determinant 
\begin{equation}\label{defn_Delta_d}
    \bm{V}_{\bm{d}}=\bm{V}_{\bm{d}}(\bm{x},\bm{y}):=\det( x_i^a y_i^b)_{\myatop{1\leq i\leq n}{(a,b)\in\bm{d}}}
 \end{equation}
of the matrix whose entries are monomials $x_i^a y_i^b$, where the indices of the variables correspond to row number, and columns are indexed by cells $(a,b)$ (elements) lying in $\bm{d}$. These may be ordered lexicographically for the purpose of the definition\footnote{Any order would work.}. Using the Reynold operator 
		$$\rho_n:=\frac{1}{n!} \sum_{\sigma\in \S_n} \sgn(\sigma)\,\sigma,$$
and ordering cells of $\bm{d}=\{(a_1,b_1),(a_2,b_2),\ldots, (a_n,b_n)\}$ as above, we may also write
   $$ \bm{V}_{\bm{d}}=\rho_n(\bm{x}^{\bm{a}}\bm{y}^{\bm{b}}),\qquad{\rm with}\qquad \bm{a}:=(a_1,a_2,\ldots, a_n)\quad {\rm and}\quad \bm{b}:=(b_1,b_2,\ldots, b_n).$$
Then one defines $\M_{\bm{d}}$ to be the linear span of all partial derivatives (to any order) of $\bm{V}_{\bm{d}}$ with respect to the variables $x_i$ and $y_j$. When $\bm{d}$ is the Ferrers diagram of a partition $\mu$, one denotes by $\M_\mu$ the resulting module.
It has been shown by Mark Haiman (see~\cite{haimanvanishing}) that
  \begin{equation}
      \MacH_\mu(q,t;\bm{z})=\M_\mu(q,t;\bm{z}).
  \end{equation}
 This is known as the $n!$-theorem. In particular, the $(k,j)$-graded component $\M_\mu^{(k,j)}$ is described as
    \begin{equation}\label{M_mu_homog_component}
       \M_\mu^{(k,j)}:=\mathcal{L}\{\partial\bm{x}^\alpha\partial\bm{y}^\beta \bm{V}_\mu\ |\ 
              \alpha,\beta\in \N^n,\ {\rm with}\  |\alpha| = \eta(\mu')-k,\ {\rm and}\  |\beta| = \eta(\mu)-j\},
    \end{equation}
 with $\mathcal{L}$ standing for ``linear span of'', and using vector notation, so that $\partial\bm{x}^\alpha=\partial x_1^{a_1}\cdots \partial x_n^{a_n}$ if $\alpha=(a_1,\ldots, a_n)$ (likewise for $\partial\bm{y}^\beta$). Observe that, for the special case $\mu=(n)$, the determinant $\bm{V}_{(n)}$ is exactly the Vandermonde determinant. The module $\M_{(n)}$ is thus the classical module of $\S_n$-harmonic polynomials (isomorphic to the $\S_n$-coinvariant module, aka the cohomology ring of the full flag manifold).
 
 \subsection*{Moving boxes} There is a nice description of the effect on $\bm{V}_{\bm{d}}$ of any symmetric operator of the form 
 	$$\kappa_{jk}:=\sum_{i=1}^n \partial y_i^j \partial x_i^k$$
for any $k+j\geq 1$. Indeed, one has
\begin{equation}
   \kappa_{jk}( \bm{V}_{\bm{d}} )= \sum_{(a,b)\in \bm{d}}  \gamma_{ab}^{jk}\, \bm{V}_{{\bm{d}_{ab}^{jk}}},
\end{equation}
where ${\bm{d}_{ab}^{jk}}$ is the diagram obtained from $\bm{d}$ by replacing the cell $(a,b)$ by the cell $(a-j,b-k)$, with\footnote{Setting $(a)_j:=a\,(1-1)\cdots (a-j+1)$.}
  $$\gamma_{ab}^{jk}:=\begin{cases}
     0,& \text{if}\ (a-j,b+k)\notin \N\times \N, \\
     0, & \text{if}\ (a-j,b+k)\in \bm{d},\\
     \pm (a)_{j} (b)_{k} & \text{otherwise}.
\end{cases}$$
In this last case, the sign corresponds to the change of the order of cells that occurs when one replaces $(a,b)$ by $(a-j,b-k)$.
This may easily be shown using the fact that the Reynold operator commutes with $\kappa_{jk}$.
Using this interpretation, observe that for all partition $\mu$ one has $\kappa_{jk} \bm{V}_\mu=0$. In a sense, this says that $ \bm{V}_\mu$ is a diagonal harmonic polynomial. The operator $\kappa_x:=\kappa_{10}$, simply ``slides'' cells to the left (up to some constant factor). 

\subsection*{$W$-expansions of $q^{\alpha}Q_{mn}(q,1/q;\bm{z})$}~
\bigskip
\begin{center}
\hrule
  {$n=3$}
 \end{center}
{\small $\begin{aligned} 
	& m=1:\qquad W_{111},\\
	& m=2:\qquad  W_{111} + qW_{21},\\
	& m=3:\qquad  [3]_qW_{111} + q\,[2]_q[3]_qW_{21} + q^{3}W_{3},\\
	& m=4:\qquad  W_{111} + q\,[3]_qW_{21} + q^{3}W_{3};
\end{aligned}$}
\begin{center}
\hrule
  {$n=4$}
 \end{center}
{\small  $\begin{aligned} 
	& m=1:\qquad  W_{1111},\\
	& m=2:\qquad  [2]_qW_{1111} + q\,[2]_qW_{211} + q^{2}W_{22},\\
	& m=3:\qquad  W_{1111} + q\,[2]_qW_{211} + q^{2}W_{22} + q^{3}W_{31},\\
	& m=4,\qquad  [4]_qW_{1111} + q\,[3]_q[4]_qW_{211} + q^{2}\tqbinom{4}{2}W_{22} 
		+ q^3\,[3]_q[4]_qW_{31} + q^{6}W_{4},\\
	& m=5,\qquad W_{1111} + q\, [4]_qW_{211} + q^{2} ([4]_q/{ [2]_q})W_{22} 
		+ q^{3}\tqbinom{4}{2}W_{31} + q^{6}W_{4};
\end{aligned}$}
\smallskip
\hrule
\begin{center}
  {$n=5$}
 \end{center}
{\small $\begin{aligned} 
	& m=1:\qquad  W_{11111},\\
	& m=2:\qquad  W_{11111} + qW_{2111} + q^{2}W_{221},\\
	& m=3:\qquad  W_{11111} + q\,[2]_qW_{2111} + q^2\,[2]_qW_{221} + q^{3}W_{311} + q^{4}W_{32},\\
	& m=4:\qquad  W_{11111} +q\,[3]_qW_{2111} + q^2\,[3]_qW_{221} + q^3\,[3]_qW_{311} 
		+ q^4\,[3]_qW_{32} + q^{6}W_{41},\\
	& m=5:\qquad [5]_q W_{11111} +q\,[4]_q[5]_q W_{2111} + q^2\,[4]_q[5]_q W_{221} 
		+ q^3[5]_q\tqbinom{4}{2}W_{311} \\
		&\hskip6cm+ q^4[5]_q\tqbinom{4}{2}W_{32} + q^{6}[4]_q [5]_q W_{41} +q^{10} W_5,\\
	& m=6:\qquad W_{11111} + q\,[5]_qW_{2111} + q^2\,[5]_qW_{221} 
		+q^{3}\tqbinom{5}{2}W_{311} 
		+ q^{4}\tqbinom{5}{2}W_{32} 
		+ q^{6}\tqbinom{5}{2}W_{41} 
		+ q^{10}W_{5}.
\end{aligned}$}
\smallskip
\hrule
\begin{center}
  {$n=6$}
 \end{center}
{\small $\begin{aligned} 
	& m=1:\qquad  W_{111111},\\
	& m=2:\qquad  [2]_qW_{111111} + q\,[2]_qW_{21111} +q^2[2]_qW_{2211} + q^{3}W_{222},\\
	& m=3:\qquad  [3]_qW_{111111} + q\,[2]_q[3]_qW_{21111}  + q^2\,[2]_q[3]_qW_{2211} 
		+ q^3[3]_qW_{222}  + q^3[3]_qW_{3111}  \\
		&\hskip6cm +q^4\,[2]_q[3]_qW_{321}  + q^{6}W_{33},\\
	& m=4:\qquad  [2]_qW_{111111}  + q\,[2]_q[3]_qW_{21111} +q^2\,[2]_q[3]_qW_{2211} 
		+ q^3[3]_qW_{222}+ q^3\,[2]_q[3]_qW_{3111}\\
		&\hskip6cm  + q^4\,[2]_q^2[3]_qW_{321}  +q^6\,[2]_qW_{33} + q^6\,[2]_qW_{411} + q^7\,[3]_qW_{42},\\
	& m=5:\qquad  W_{111111} +q\,[4]_qW_{21111} +q^2[4]_qW_{2211} + q^3(\tqbinom{4}{2}/[3]_q) W_{222} 
	+ q^3\tqbinom{4}{2}W_{3111} \\
		&\hskip6cm   + q^4[3]_q[4]_qW_{321}+ q^6([4]_q/ [2]_q)W_{33} +q^6[4]_qW_{411} \\
		&\hskip6cm   + q^7\tqbinom{4}{2}W_{42} + q^{10}W_{51},\\
	& m=6:\qquad [6]_qW_{111111} + q\,[5]_q[6]_qW_{21111} + q^2[5]_q[6]_qW_{2211} +q^2\tqbinom{6}{2}W_{222}\\
		&\hskip6cm   + q^3[3]_q\tqbinom{6}{2}W_{3111} +  q^4([4]_q[5]_q[6]_q)W_{321} \\
		&\hskip6cm    + q^6\tqbinom{6}{3}W_{33} + q^7([3]_q[4]_q[5]_q[6]_q/[2]_q^2)W_{42}\\
		&\hskip6cm     + q^{10}[5]_q[6]_qW_{51} + q^{15}W_{6},\\
	& m=7:\qquad  W_{111111} +q\, [6]_qW_{21111} + q^2\, [6]_qW_{2211} + q^3[3]_qW_{222} 
			+ q^3\textstyle \qbinom{6}{2}W_{3111} + q^4([5]_q[6]_q)W_{321}\\
		&\hskip6cm     +q^6(\tqbinom{6}{3}/[4]_q) W_{33} 
		+ q^6\tqbinom{6}{3}W_{411}\\
		&\hskip6cm     + q^6([4]_q[5]_q[6]_q/[2]_q^2)W_{42} +q^{10}\tqbinom{6}{2}W_{51} + q^{15}W_{6}.
\end{aligned}$}
\smallskip
\hrule

%%%%%%%%%%%%%%%%%%%%%%%%%%%%%%%%%%%%%%%%%%%%
\renewcommand{\refname}{\bleu{References}}

\end{document}